\documentclass[11pt,oneside,dvipsnames]{amsart}
\usepackage{stackrel}
\usepackage[french, english]{babel}
\usepackage{a4wide}
\usepackage[utf8,latin1]{inputenc}
\usepackage{times}
\usepackage{adjustbox}
\usepackage{float}
\usepackage[cyr]{aeguill}
\usepackage{amsmath,amscd}
\usepackage{amssymb}
\usepackage{mathrsfs}
\usepackage{amsthm}
\usepackage{geometry}
\usepackage{graphicx}
\usepackage{epstopdf}
\usepackage{amsfonts}
\usepackage{amsmath}
\usepackage{amsmath,mathtools}
\usepackage{amsmath,graphicx}
\usepackage{amscd}
\usepackage{dsfont}
\usepackage{latexsym}
\usepackage{verbatim}
\usepackage{comment}
\usepackage{tikz-cd}
\usetikzlibrary{matrix,arrows}
\usepackage{tikz}
\usetikzlibrary{tikzmark}
\usepackage{tensor}
\usepackage{leftidx}
\usepackage[retainorgcmds]{IEEEtrantools}
\usepackage[title]{appendix}
\usepackage[mathscr]{euscript}
\usepackage[colorinlistoftodos]{todonotes}
\usepackage{bbm}
\usepackage{url}
\usepackage{upgreek}
\usepackage{mathtools, nccmath}
\usepackage[utf8]{inputenc}
\usepackage{multirow}
\usepackage{rotating}
\usepackage{marginnote}
\usepackage{hyperref} 
\usepackage{nameref}
\usepackage{ytableau}
\usepackage{youngtab}

\usepackage[all]{xy}
\usepackage[T1]{fontenc}

\usepackage{hyperref}
\usepackage{array}
\usepackage{tabularx}
\usepackage{yfonts}
\usepackage{setspace}  
 \usepackage[foot]{amsaddr}
\makeatletter
\renewcommand\subsubsection{\@startsection{subsubsection}{3}{\z@}%
                                     {-2ex\@plus -1ex \@minus -.2ex}%
                                     {-0.5em}
                                     {\normalfont\normalsize\bfseries}}
                                     
                                     \renewcommand\subsection{\@startsection{subsection}{3}{\z@}%
                                     {-3.25ex\@plus -1ex \@minus -.2ex}%
                                     {-0.5em}
                                     {\normalfont\normalsize\bfseries}}

\makeatother

\theoremstyle{plain}
\newtheorem{Thm}{Theorem}
\newtheorem*{Thm*}{Theorem}
\newtheorem{ThmA}{Theorem}

\newtheorem{Prop}[Thm]{Proposition}
\newtheorem{Lem}[Thm]{Lemma}
\newtheorem{Cor}[Thm]{Corollary}

\newtheorem*{Con*}{Conjecture}

\theoremstyle{definition}

\newtheorem{Rem}[Thm]{Remark}

\newtheoremstyle{named}{}{}{}{}{\bfseries}{.}{.5em}{\thmnote{#3}#1}
\theoremstyle{named}

\DeclareMathOperator{\Coeff}{Coeff}

\DeclareMathOperator{\End}{End}
\DeclareMathOperator{\Exp}{Exp}

\DeclareMathOperator{\gl}{\mathfrak{gl}}
\DeclareMathOperator{\GL}{GL}

\DeclareMathOperator{\Gr}{Gr}
\DeclareMathOperator{\Hom}{Hom}

\DeclareMathOperator{\Id}{Id}
\DeclareMathOperator{\Ima}{Im}

\DeclareMathOperator{\lisom}{\!\!\smash{\begin{array}{c}\sim\\[-1em]
\longrightarrow\end{array}}\!\!}
\DeclareMathOperator{\lisomL}{\!\!\smash{\begin{array}{c}\sim\\[-1em]
\longleftarrow\end{array}}\!\!}

\DeclareMathOperator{\Ker}{Ker}

\DeclareMathOperator{\pr}{pr}

\DeclareMathOperator{\rk}{rk}

\DeclareMathOperator{\Rep}{Rep}

\DeclareMathOperator{\ds}{\!/\mkern-2mu/\mkern-2mu}

\newcommand{\hooklongrightarrow}{\lhook\joinrel\longrightarrow}

\makeatletter
\newcommand*{\rom}[1]{\expandafter\@slowromancap\romannumeral #1@}
\makeatother

\renewcommand{\descriptionlabel}[1]{\hspace\labelsep\upshape\bfseries #1.}
\makeatletter
\let\orgdescriptionlabel\descriptionlabel
\renewcommand*{\descriptionlabel}[1]{%
  \let\orglabel\label
  \let\label\@gobble
  \phantomsection
  \edef\@currentlabel{#1}%
  \let\label\orglabel
  \orgdescriptionlabel{#1}%
}
\makeatother

\usepackage{subfig}	 

\title{The irreducibility of tame character varieties}
\author{Cheng Shu}
\address[C. Shu]{School of Mathematical Sciences, Tongji University, Shanghai, China}
\email[C. Shu]{shuc@tongji.edu.cn}

\usepackage{pxfonts}
\usepackage{indentfirst}
\usepackage{setspace}  
\usepackage{hyperref}
\colorlet{ivory}{Apricot!30!}
\colorlet{space}{black!85!}
\definecolor{bgc}{RGB}{29, 44, 46}
\definecolor{txt}{RGB}{223, 222, 189}
\definecolor{cmd}{RGB}{206, 151, 88}
\begin{document}
\let\bs\boldsymbol
\begin{abstract}
Tame character varieties are known to be irreducible when the parametrised local systems have monodromies lying in closures of generic  conjugacy classes; we will establish their irreducibility without the genericity assumption. A combination of nonabelian Hodge correspondence and variation of parabolic weights reduces the problem to the case of  almost generic semi-simple monodromies, which is then treated using a point-counting result of Scognamiglio for nongeneric character stacks. The connectedness result is then used to clarify a couple of facts concerning the symplectic leaves.
\end{abstract}
\maketitle

\tableofcontents
\addtocontents{toc}{\protect\setcounter{tocdepth}{-1}}
\setcounter{tocdepth}{1}
\numberwithin{Thm}{section}
\numberwithin{equation}{section}
\addtocontents{toc}{\protect\setcounter{tocdepth}{1}}
\section{Introduction}\label{INT}

Character varieties are moduli spaces of local systems on Riemann surfaces and form an important family of Poisson varieties. Symplectic moduli spaces are obtained by restricting the monodromies of the local systems at the punctures to given conjugacy classes. We will call these symplectic varieties the tame character varieties, which are the object of study of this paper. The tame property amounts to saying that the connections corresponding to the parametrised local systems have regular singularities. We will not consider their wild counterpart as studied in \cite{BB} and \cite{Boa14}.

It is an open problem whether these symplectic varieties are connected or even irreducible. By counting points over finite fields, Hausel-Letellier-Rodriguez-Villegas \cite{HLRV} \cite{HLRV2}, continuing on the earlier work \cite{HRV}, obtained explicit results on the cohomology of character varieties, which show in particular that if the prescribed monodromies are generic semi-simple, then tame character varieties are connected. Subsequent work of Letellier \cite{Le15} established the connectedness without requiring the monodromies to be semi-simple, but the genericity assumption remained essential.

These results agree with the general geometric principle (like the theorem of Bertini) that a variety cut out by generic equations are smooth and connected, while in degenerate cases it can be singular or have unexpected connected components, possibly of unexpected dimension. In this sense, the main theorem of the present article is surprising since it says that nongeneric defining equations for character varieties never produce unexpected components, although the dimension may be larger than expected. It seems to be difficult to establish the connectedness in nongeneric cases using purely algebro-geometric methods. As far as the author is aware of, the connectedness of nongeneric $\GL_n$-character varieties is only known in the following cases:
\begin{itemize}
\item Arbitrary $g$ and $k$. In \cite{Ina}, Inaba proved the irreducibility of some smooth moduli spaces of stable parabolic connections, but more general singular moduli spaces were left unaddressed. Moreover, the parabolic structures are required to be full flags. See \cite[\S 5]{Ina} for the precise statements.
\item $g>0$, $k=1$, and the monodromy at the puncture is central but not necessarily generic. One proof was given by Li using the infinite dimensional construction of the moduli spaces of connections; see \cite{Li}. Another proof due to de Cataldo requires $g>1$ and uses the structure of the Hitchin map; see \cite{dC}.
\item A recent preprint of Lee-Lee \cite{LL} also proved the connectedness of the locus of irreducible local systems with a particular constraint on the eigenvalues of monodromies.
\item Nongeneric multiplicative quiver varieties for some particular quivers were considered in \cite[Section 4.4]{KS2}; they also discussed the connectedness of representation type strata, which is the topic of Section \ref{sec-Symp-L} of the present paper.
\end{itemize}

In a recent study of the (multiplicative) Deligne-Simpson problem \cite{Shu25}, it was realised that in certain cases it is possible to produce a surjective continuous map from a moduli space with generic parameters to a given moduli space with nongeneric parameters, thus reducing the connectedness problem to known cases. In this article, we will use this strategy to prove the connectedness of tame character varieties in full generality. However, the result of Hausel-Letellier-Rodriguez-Villegas is not sufficient. In general, we can only produce a surjective map from a moduli space with what we call \textit{almost} generic parameters, and the connectedness of such a space cannot be further reduced to simpler spaces. The solution turns out to be an application of a more general result of Scognamiglio \cite{Sco} on the point-counting of nongeneric character stacks. 

\subsection{The main result}\label{subsec-main-thm}\hfill

Our main result concerns the moduli spaces $\mathcal{M}_B(\mathbf{P},\bar{\mathcal{C}},\bs\beta)$ which we describe in more detail now. Let $\bar{X}$ be a compact Riemann surface of genus $g\ge 0$ and let $X=\bar{X}\setminus\{x_1,\ldots,x_k\}$ be the open Riemann surface obtained by removing $k$ points. A filtered local system on $X$ is a local system $\mathcal{L}$ equipped with a filtered structure $\mathcal{L}_{[\bullet,\bullet]}$, which consists of a flag of locally defined local systems $\mathcal{L}_{[j,\bullet]}$ for every $1\le j\le k$; for each $j$, we choose a simply connected open neighbourhood $U_j$ of $x_j$, and the flag $\mathcal{L}_{[j,\bullet]}$ is defined on $U_j^{\ast}=U_j\setminus\{x_j\}$. We may write $\mathcal{L}_{[j,\bullet]}$ more explicitly as 
$$
\mathcal{L}|_{U_j^{\ast}}=\mathcal{L}_{[j,0]}\supset\mathcal{L}_{[j,1]}\supset\cdots\supset\mathcal{L}_{[j,\mu_j+1]}=0
$$
for some $\mu_j$. The choice of $U_j$ in the definition of filtered local system will not matter. We need to fix the following data in the construction of a moduli space:
\begin{itemize}
\item The type of the filtered structure $\mathcal{L}_{[\bullet,\bullet]}$ is represented by a tuple of parabolic subgroups $\mathbf{P}=(P_j)_j$, where each $P_j\subset\GL_n$ is conjugate to the parabolic subgroup preserving the fibre of $\mathcal{L}_{[j,\bullet]}$ at some point of $U_j^{\ast}$.
\item The graded monodromy $\mathcal{C}$. Choose a point $y_j\in U_j^{\ast}$, then the monodromy of $\mathcal{L}|_{U_j^{\ast}}$ is an element $M_j\in\GL(\mathcal{L}_{y_j})$, where $\mathcal{L}_{y_j}$ is the fibre of $\mathcal{L}$ over $y_j$. Since $M_j$ preserves the flag of vector spaces $\mathcal{L}_{[j,\bullet],y_j}$, it induces an element $\tilde{M}_j\in\prod_{i=0}^{\mu_j}\GL(\mathcal{L}_{[j,i],y_j}/\mathcal{L}_{[j,i+1],y_j})$. We will call $\tilde{M}_j$ the graded monodromy of $\mathcal{L}$ at $x_j$; it is required to lie in the Zariski closure $\bar{C}_j$ of a given conjugacy class $C_j$ of $\prod_{i=0}^{\mu_j}\GL(\mathcal{L}_{[j,i],y_j}/\mathcal{L}_{[j,i+1],y_j})$. The tuple of conjugacy classes $(C_j)_{1\le j\le k}$ will be denoted by $\mathcal{C}$ and the tuple of their closures by $\bar{\mathcal{C}}$.
\item The stability condition $\bs\beta$. For each $j$, let $(\beta_{[j,i]})_{0\le i\le \mu_j}$ be a sequence of strictly increasing real numbers and write $\bs\beta=(\beta_{[j,i]})_{j,i}$. The $\bs\beta$-degree of a filtered local system $(\mathcal{L},\mathcal{L}_{[\bullet,\bullet]})$ is defined to be
$$
\deg_{\bs\beta}\mathcal{L}:=\sum_{1\le j\le k}\sum_{0\le i\le \mu_j}\beta_{[j,i]}(\rk\mathcal{L}_{[j,i]}-\rk\mathcal{L}_{[j,i+1]}).
$$ 
We say that a filtered local system $(\mathcal{L},\mathcal{L}_{[\bullet,\bullet]})$ of $\bs\beta$-degree zero is $\bs\beta$-semi-stable (resp. stable) if for any local system $\mathcal{M}\subset\mathcal{L}$ with its induced filtered structure $\mathcal{M}_{[\bullet,\bullet]}$ and induced weight $\bs\beta$, we have
\begin{equation}\label{eq-filt-locsys-stability}
\deg_{\bs\beta}\mathcal{M}\le0\text{ (resp. $<0$).}
\end{equation}
\end{itemize}
Recall that for any $j$, the flag $\mathcal{M}_{[j,\bullet]}$ consists of the distinct local systems on $U_j^{\ast}$ among $\mathcal{L}_{[j,\bullet]}\cap \mathcal{M}|_{U_j^{\ast}}$, and the local subsystem $\mathcal{M}_{[j,i']}$ is assigned the weight $\beta_{[j,i]}$ if $\mathcal{L}_{[j,i]}\cap\mathcal{M}|_{U_j^{\ast}}=\mathcal{M}_{[j,i']}$ and $\mathcal{L}_{[j,i+1]}\cap\mathcal{M}|_{U_j^{\ast}}\ne\mathcal{M}_{[j,i']}$.

There exists a coarse moduli space $\mathcal{M}_B(\mathbf{P},\bar{\mathcal{C}},\bs\beta)$ of $\bs\beta$-semi-stable filtered local systems on $X$ of zero $\bs\beta$-degree and of type $\mathbf{P}$ with graded monodromies lying in $\bar{\mathcal{C}}$. The main theorem below does not impose any restriction on the invariants $\mathbf{P}$, $\bar{\mathcal{C}}$ and $\bs\beta$ other than what is necessary for the moduli space to be nonempty.

\begin{ThmA}\label{Thm-main}
Suppose that $\mathcal{M}_B(\mathbf{P},\bar{\mathcal{C}},\bs\beta)$ is nonempty. Then, it is an irreducible normal algebraic variety.
\end{ThmA}
\begin{Rem}
Indeed, it follows from our proof and \cite[Theorem 5.4]{KS} that $\mathcal{M}_B(\mathbf{P},\bar{\mathcal{C}},\bs\beta)$ has symplectic singularities.
\end{Rem}
\begin{Rem}
It is essential to take closures of conjugacy classes in the definition of moduli spaces in order for the irreducibility to hold. The computation of Kamgarpour-Nam-Pusk\'as showed that if we impose unipotent monodromies without taking closures, then we may end up with a disconnected character variety; see \cite[Section 7.1.3]{KNP}. Adding boundary strata allows different orbits to contain a common closed orbit in their respective closures. Indeed, the two local systems in \textit{loc. cit.} representing the two connected components have the same semi-simplification, and thus they represent the same point in our setting. Beware that nontrivial unipotent conjugacy classes are not affine, and thus taking the quotient (even for a free action) of a representation variety defined by such conjugacy classes in general only produces an algebraic space instead of a scheme; of course, this does not affect the discussion on connected components in \textit{op. cit.}.
\end{Rem}

\subsection{The method of proof}\hfill

Our proof combines two perspectives of the same moduli space. As a moduli space of filtered local systems, the space $\mathcal{M}_B(\mathbf{P},\bar{\mathcal{C}},\bs\beta)$ fits in the nonabelian Hodge theory, and thus transcendental maps can be used to relate moduli spaces defined by different monodromies. However, we may also regard $\mathcal{M}_B(\mathbf{P},\bar{\mathcal{C}},\bs\beta)$ as a multiplicative quiver variety $\mathcal{M}(\bs\zeta,\mathbf{d},\bs\theta)$ (or an open subset thereof if the genus of $X$ is not zero). On the one hand, this allows us to apply the results of Kaplan-Schedler to deduce that it is a normal variety. On the other hand, there are reflection isomorphisms from the theory of quiver representations, relating moduli spaces defined for different dimension vectors. The latter feature reduces the connectedness problem to those dimension vectors that are imaginary roots in the Kac-Moody root system associated to the underlying quiver.

If $\mathbf{d}$ is an imaginary root that is indivisible relative to the eigenvalues, then the method of \cite[Proposition 4.2]{Shu25} can be used to show that such a character variety is connected. The proof, as mentioned above, goes in two steps:
\begin{itemize}
\item Construct a continuous (transcendental) surjective map 
\begin{equation}\label{eq-intro-surj}
\mathcal{M}\rightarrow\mathcal{M}_B(\mathbf{P},\bar{\mathcal{C}},\bs\beta),
\end{equation} 
where $\mathcal{M}$ is a character varieties defined by generic semi-simple monodromies.
\item Apply the theorem of Hausel-Letellier-Rodriguez-Villegas that we mentioned above to deduce the connectedness of $\mathcal{M}$.
\end{itemize}
However, for a multiple $m\mathbf{d}$ of such a dimension vector and the same eigenvalues, this method fails and essential difficulties arise. The property of $m\mathbf{d}$ that is crucial for the discussion of the connectedness is whether it is isotropic or anisotropic; see \S \ref{subsec-KM-roots} for the definitions.

If $\mathbf{d}$ is isotropic, then the main problem lies in proving a map like (\ref{eq-intro-surj}) to be surjective. In fact, the stable locus of $\mathcal{M}_B(\mathbf{P},\bar{\mathcal{C}},\bs\beta)$ turns out to be empty, and thus it is unlikely that variation of stability condition defines a surjective map. However, the connectedness of $\mathcal{M}_B(\mathbf{P},\bar{\mathcal{C}},\bs\beta)$ does hold according to \cite[Theorem 4.4]{Shu25}. Indeed, this theorem shows that it is the $m$-th symmetric product of a character variety defined by a relatively indivisible dimension vector, which implies the connectedness. 

The focus of this article is the anisotropic dimension vectors. What continues to be working is the nonabelian Hodge theory and variation of stability techniques that reduce the connectedness problem to character varieties with semi-simple monodromies, resulting in a surjective map like (\ref{eq-intro-surj}). The new phenomenon in this case is that $\mathcal{M}$ is no longer defined by generic monodromies because of the divisibility of $\mathbf{d}$. Moreover, we cannot reduce the connectedness problem to simpler spaces, so that it is as close to a generic character variety as possible. Fortunately, the point-counting techniques can be suitably generalised to tackle this situation, thanks to previous work of Scognamiglio. It can be read from his point-counting formula that the number of $\mathbb{F}_q$-points of the corresponding representation variety is a monic polynomial in $q$, which gives the expected connectedness. We remark that Scognamiglio's formula does apply to character varieties defined by arbitrary nongeneric semi-simple conjugacy classes. However, a direct computation of the leading term may lead to complicated (perhaps impossible) combinatorics, and the use of nonabelian Hodge theory reduces our computational load to a minimum.

\subsection{Relation to other work}\hfill

Our main result addresses the connectedness problem for a particular family of multiplicative quiver varieties. For the quivers of affine type $\tilde{A}_l$ with an extending vertex, the connectedness of the associated multiplicative quiver varieties with generic deformation parameters and the trivial stability condition was proved by Braverman-Etingof-Finkelberg (see \cite[\S 5.3]{BEF}). Of course, it is natural to ask if all multiplicative quiver varieties are connected. Indeed, this was listed as an open problem in \cite[\S 7]{ST}. (The authors gave a short explanation for the irreducibility in the comet-shaped quiver case, which is based on \cite{Ina}; however, as we remarked above, Inaba's results assume the parabolic structures to be full flags, and it is not clear how his arguments can be generalised.) The additive version of this problem with the trivial stability condition was solved by Crawley-Boevey in \cite{CB01} and \cite{CB02}.  Nakajima studied the connectedness of Nakajima quiver varieties in \cite[\S 6]{Naka98} and \cite[\S 7.5]{Naka01}, which turned out to be an implication of Crawley-Boevey's results. Using hyperk\"ahler twisting, Bellamy-Schedler deduce from Crawley-Boevey's results that quiver varieties with nontrivial stability conditions are also connected; see \cite[Proposition 3.21]{BS0}. This in turn allowed them to show that the representation type strata of quiver varieties are irreducible and are precisely the symplectic leaves (the loci where the induced Poisson structure is nondegenerate); see \cite[Corollary 3.25]{BS0}. We will prove an analogous result for tame character varieties in the last section.

\subsection*{AI disclosure}\hfill

Various versions of ChatGPT were used in the preparation of this article, including 5.5 Pro, 5.6 Pro and 6 Astra Pro. The chatbot was used to collect and clarify a couple of facts concerning the quiver varieties for quivers with framing vertices and to find some references which are included in this version of the article. It also provided Remark \ref{Rem-CaEs-paper}. This article contains no AI-generated texts or ideas. We
take full responsibility for the final content of this work.

\subsection*{Acknowledgement}\hfill

This article was initiated while I was visiting Xueqing Wen at Chongqing University of Technology in October 2025, and was completed while visiting him again in November 2026. I would like to thank him for his hospitality during my stays. I also thank Emmanuel Letellier for his interest and a couple discussions. This work is supported by the China Postdoctoral Science Foundation under Grant Numbers 2025M773090 and 2025T180845.

\section{Character varieties and multiplicative quiver varieties}\label{sec-Char-Quiv}

In this section, we review the multiplicative quiver varieties introduced by Crawley-Boevey-Shaw \cite{CBS} as well as the tame character varieties (i.e., moduli spaces of the filtered local systems introduced by Simpson \cite{Si90}). The relation between these varieties has been clarified by Yamakawa \cite{Yama}, Schedler-Tirelli \cite{ST} and Scognamiglio \cite{Sco}, and we know that tame character varieties can be identified with open subvarieties of multiplicative quiver varieties. However, each of these papers imposes some constraints on the parametrised local systems. Yamakawa's results are restricted to the genus zero case, and the graded monodromies are required to be central. Results of Schedler-Tirelli and Scognamiglio are only stated for affine character varieties, so that trivial filtered structures are put on local systems. Of course, the proof for the general case is completely analogous to these existing results. We prefer giving some more details to making it an assertion. Besides, a general form of reflection isomorphisms between multiplicative quiver varieties was stated in \cite[\S 2.5]{ST}, and we also take this opportunity to provide more details of its proof. In fact, it is not completely free of effort to see that Yamakawa's proof for the stable loci works in the generality as stated there. I hope that what we provide here could convince the readers that reflection isomorphisms, defined for the ambient multiplicative quiver varieties, do preserve the tame character varieties when $g>0$.

\subsection{Multiplicative quiver varieties}\hfill

A quiver $Q=(Q_0,Q_1)$ consists of a set of vertices $Q_0$ and a set of arrows $Q_1$. For any $a\in Q_1$, we denote by $t(a)\in Q_0$ the tail of $a$ and by $h(a)\in Q_0$ the head of $a$. A representation of $Q$ is a $Q_0$-graded vector space $\mathbf{V}=(V_v)_{v\in Q_0}$ together with a collection of linear maps $\bs\phi=(\phi_a)_{a\in Q_1}$ where $\phi_a\in\Hom(V_{t(a)},V_{h(a)})$; a representation will be denoted by $(\mathbf{V},\bs\phi)$, and either component of this pair may be omitted from the notation if the context permits. Let $Q_1^{\ast}$ be a set of arrows which contains for each $a\in Q_1$ a unique arrow $a^{\ast}$ going from $h(a)$ to $t(a)$. We will write $\bar{Q}_1=Q_1\sqcup Q_1^{\ast}$ and the double of $Q$ is the quiver $\bar{Q}=(Q_0,\bar{Q}_1)$. For any $\mathbf{d}=(d_v)_{v\in Q_0}\in(\mathbb{Z}_{\ge0})^{Q_0}$, the space of representations of $\bar{Q}$ of dimension vector $\mathbf{d}$ is the vector space
$$
\Rep(\bar{Q},\mathbf{d}):=\bigoplus_{a\in \bar{Q}_1}\Hom(\mathbb{C}^{d_{t(a)}},\mathbb{C}^{d_{h(a)}}).
$$
We say that a quiver representation $\bs\phi=(\phi_a)_{a\in\bar{Q}_1}\in\Rep(\bar{Q},\mathbf{d})$ is invertible if for every $a\in Q_1$, we have $\det(\Id+\phi_a\phi_{a^{\ast}})\ne0$. Note that $\det(\Id+\phi_a\phi_{a^{\ast}})\ne0$ is equivalent to $\det(\Id+\phi_{a^{\ast}}\phi_a)\ne0$. We will denote by 
$$
\Rep^{\circ}(\bar{Q},\mathbf{d})\subset\Rep(\bar{Q},\mathbf{d})
$$
the open subvariety consisting of invertible representations. Let $\epsilon:Q_1\sqcup Q_1^{\ast}\rightarrow\{\pm1\}$ be the map defined by $\epsilon(a)=1$ if $a\in Q_1$ and $\epsilon(a)=-1$ otherwise. Choose a total ordering $\prec$ on $\bar{Q}_1$ and define
\begingroup
\allowdisplaybreaks
\begin{align*}
\mu:\Rep^{\circ}(\bar{Q},\mathbf{d})&\longrightarrow\prod_{v\in Q_0}\GL_{d_v}=:G\\
\bs\phi&\longmapsto (\prod_{\substack{a\in \bar{Q}_1\\h(a)=v}}^{\prec}(\Id+\phi_a\phi_{a^{\ast}})^{\epsilon(a)})_{v\in Q_0},
\end{align*}
\endgroup
where the factors of the product $\prod^{\prec}$ are placed in the order determined by the ordering $\prec$; this is the quasi-Hamiltonian moment map for $\Rep^{\circ}(\bar{Q},\mathbf{d})$ (see \cite[\S 3.1]{Yama}). The component of $\mu$ corresponding to $v$ will be denoted by $\mu_v$.

Let $\bs\phi\in\Rep(\bar{Q},\mathbf{d})$ be a representation. A subrepresentation of $\bs\phi$ is a $Q_0$-graded vector space $(V_v)_{v\in Q_0}$ with $V_v\subset\mathbb{C}^{d_v}$ and $\phi_a(V_{t(a)})\subset V_{h(a)}$ for any $a\in\bar{Q}_1$.
Let $\bs\theta=(\theta_v)_{v\in Q_0}\in\mathbb{R}^{Q_0}$ be a vector satisfying 
$$
\bs\theta\cdot\mathbf{d}:=\sum_{v\in Q_0}\theta_vd_v=0.
$$
We say that $\bs\phi$ is $\bs\theta$-stable (resp. semi-stable) if for any nontrivial proper subrepresentation $(V_v)_{v\in Q_0}$ we have
$$
\sum_{v\in Q_0}\theta_v\dim V_v<0 \text{ (resp. $\le0$)}.
$$
The open subset of $\bs\theta$-semi-stable invertible representations will be denoted by 
$$
\Rep^{\circ}(\bar{Q},\mathbf{d})^{\bs\theta-ss},
$$
and for any $\bs\zeta=(\zeta_v)_{v\in Q_0}\in(\mathbb{C}^{\ast})^{Q_0}$, regarded as a central element of $G$, the $\bs\theta$-semi-stable part of the fibre $\mu^{-1}(\bs\zeta)$ will be denoted by
$$
\Rep^{\circ}(\bar{Q},\bs\zeta,\mathbf{d})^{\bs\theta-ss}.
$$
Consider the action of $G$ on $\Rep(\bar{Q},\mathbf{d})$ defined by
$$
(\phi_a)_{a\in\bar{Q}_1}\longmapsto(g_{h(a)}\phi_ag_{t(a)}^{-1})_{a\in \bar{Q}_1}
$$
for any $(g_v)_{v\in Q_0}\in G$. This action preserves the locus of $\bs\theta$-semi-stable representations, the locus of invertible representations, and the fibre $\mu^{-1}(\bs\zeta)$ for any $\bs\zeta$.

Finally, the multiplicative quiver variety associated to $\bs\zeta$, $\mathbf{d}$ and $\bs\theta$ is the good quotient (i.e. the quotient parametrising the closed orbits)
$$
\mathcal{M}(Q,\bs\zeta,\mathbf{d},\bs\theta):=\Rep^{\circ}(\bar{Q},\bs\zeta,\mathbf{d})^{\bs\theta-ss}\ds G.
$$
We may omit $Q$ or $\bar{Q}$ from the notation if no confusion arises.

\subsection{Kac-Moody roots and reflections}\label{subsec-KM-roots}\hfill

Given a quiver $Q$, the Euler form on $\mathbb{Q}^{Q_0}$ is defined by
$$
(\mathbf{c},\mathbf{d}):=\sum_{v\in Q_0}2c_vd_v-\sum_{a\in Q_1}c_{t(a)}d_{h(a)}-\sum_{a\in Q_1}c_{h(a)}d_{t(a)}.
$$
The coordinate vectors are written as $\mathbf{e}_v=(e_w)_{w\in Q_0}$ with $e_v=1$ and $e_w=0$ for $w\ne v$; if $v$ is a loop-free vertex, then $(\mathbf{e}_v,\mathbf{e}_v)=2$. The support of $\mathbf{d}\in(\mathbb{Z}_{\ge0})^{Q_0}$ is the subquiver obtained by removing all vertices $v$ with $d_v=0$ and the arrows connecting to these vertices. The fundamental region $\mathcal{F}(Q)\subset(\mathbb{Z}_{\ge0})^{Q_0}$ is the subset of nonzero vectors $\mathbf{d}$ having a connected support with $(\mathbf{d},\mathbf{e}_v)\le0$ for all $v$.

Let $v\in Q_0$ be a loop-free vertex. Define the following (reflection) maps
\begingroup
\allowdisplaybreaks
\begin{align*}
s_v:\mathbb{Z}^{Q_0}&\longrightarrow\mathbb{Z}^{Q_0},\\
u_v:(\mathbb{C}^{\ast})^{Q_0}&\longrightarrow(\mathbb{C}^{\ast})^{Q_0},\text{ and}\\
r_v:\mathbb{R}^{Q_0}&\longrightarrow\mathbb{R}^{Q_0},
\end{align*}
\endgroup
by
\begingroup
\allowdisplaybreaks
\begin{align*}
s_v(\mathbf{d})=\mathbf{d}',\quad&\mathbf{d}'=\mathbf{d}-(\mathbf{d},\mathbf{e}_v)e_v,\\
u_v(\bs\zeta)=\bs\zeta',\quad& \zeta'_w=\zeta_w\zeta_v^{-(\mathbf{e}_v,\mathbf{e}_w)}, \text{ and}\\
r_v(\bs\theta)=\bs\theta',\quad&\theta'_w=\theta_w-(\mathbf{e}_v,\mathbf{e}_w)\theta_v.
\end{align*}
\endgroup

The real roots are the vectors that are obtained by applying a sequence of reflections to the vectors $\pm\mathbf{e}_v$ supported on loop-free vertices; the imaginary roots are those obtained by applying a sequence of reflections to the elements of $\mathcal{F}(Q)$ up to a sign. We will denote by $R$ the set of all roots and $R^+=R\cap (\mathbb{Z}_{\ge0})^{Q_0}$ the set of positive roots. For any $\mathbf{d}$, write
$$
p(\mathbf{d})=1-\frac{1}{2}(\mathbf{d},\mathbf{d}).
$$
We have $p(\mathbf{d})=0$ if $\mathbf{d}$ is a real root. We call an imaginary root $\mathbf{d}$ an isotropic root if $p(\mathbf{d})=1$ and an anisotropic root if $p(\mathbf{d})>1$.

We say that a vector $\mathbf{d}\in(\mathbb{Z}_{\ge0})^{Q_0}$ is indivisible if it is not of the form $m\mathbf{d}_0$ for some $\mathbf{d}_0\in(\mathbb{Z}_{\ge0})^{Q_0}$ and $m\in\mathbb{Z}_{>1}$. Let $\mathbf{d}=m\mathbf{d}_0$ for some $m\in\mathbb{Z}_{>0}$ and indivisible $\mathbf{d}_0\in(\mathbb{Z}_{\ge0})^{Q_0}$. Let $\bs\zeta\in(\mathbb{C}^{\ast})^{Q_0}$ and suppose that
$$
\bs\zeta^{\mathbf{d}}:=\prod_{v\in Q_0}\zeta_v^{d_v}=1.
$$ 
We say that $\mathbf{d}$ is $\bs\zeta$-indivisible if there is no integer $0<l<m$ such that $\bs\zeta^{l\mathbf{d}_0}=1$. For a given $\mathbf{d}$ as above, a deformation parameter $\bs\zeta\in(\mathbb{C}^{\ast})^{Q_0}$ with $\bs\zeta^{\mathbf{d}}=1$ is called almost generic if every $\mathbf{d}'\in(\mathbb{Z}_{\ge0})^{Q_0}$ satisfying $\mathbf{d}'\le\mathbf{d}$ (i.e. $d'_v\le d_v$ for every $v\in Q_0$) and $\bs\zeta^{\mathbf{d}'}=1$ is a $\mathbb{Q}$-multiple of $\mathbf{d}$; if such a $\mathbf{d}'$ is either zero or $\mathbf{d}$ itself, then $\bs\zeta$ is called generic. A stability condition $\bs\theta\in\mathbb{R}^{Q_0}$ with $\bs\theta\cdot\mathbf{d}=0$ is called almost generic if every $\mathbf{d}'\in(\mathbb{Z}_{\ge0})^{Q_0}$ satisfying $\mathbf{d}'\le\mathbf{d}$ and $\bs\theta\cdot\mathbf{d}'=0$ is a $\mathbb{Q}$-multiple of $\mathbf{d}$; if such a $\mathbf{d}'$ is either zero or $\mathbf{d}$ itself, then $\bs\theta$ is called generic.

A particular set of positive roots that are of particular importance is the set $\Sigma_{\bs\zeta,\bs\theta}$ defined below. We define $R^+_{\bs\zeta,\bs\theta}:=\{\mathbf{d}\in R^+\mid \bs\zeta^{\mathbf{d}}=1,~\bs\theta\cdot\mathbf{d}=0\}$, and
\begingroup
\allowdisplaybreaks
\begin{align*}
\Sigma_{\bs\zeta,\bs\theta}:=\{\mathbf{d}\in R^+_{\bs\zeta,\bs\theta}\mid&\text{if $\mathbf{d}=\sum_{s=1}^r\mathbf{d}^{(s)}$ with $r\ge2$ and each $\mathbf{d}^{(s)}\in R^+_{\bs\zeta,\bs\theta}$,}\\
&\text{then $p(\mathbf{d})>\sum_{s=1}^rp(\mathbf{d}^{(s)})$}\}.
\end{align*}
\endgroup
We will omit $\bs\zeta$ (resp. $\bs\theta$) from the subscripts of $\Sigma_{\bs\zeta,\bs\theta}$ and $R^+_{\bs\zeta,\bs\theta}$ if $\bs\zeta=1$ (resp. $\bs\theta=0$). This set was introduced by Crawley-Boevey in his study of the Deligne-Simpson problem.

\subsection{Reflection isomorphisms}\label{subsec-Ref}\hfill

We now review the reflection isomorphisms between multiplicative quiver varieties, following Crawley-Boevey-Shaw \cite{CBS} and Yamakawa \cite{Yama}. 

Let $\mu$ and $\mu'$ be the quasi-Hamiltonian moment maps defined on $\Rep^{\circ}(\bar{Q},\mathbf{d})$ and $\Rep^{\circ}(\bar{Q},\mathbf{d}')$ respectively.
The reflection functor at a loop-free vertex $v$ sends a $\bs\theta$-semi-stable $\mathbf{d}$-dimensional invertible representation $(\mathbf{V},\bs\phi)$ satisfying $\mu(\bs\phi)=\bs\zeta$ to a $\bs\theta'$-semi-stable $\mathbf{d}'$-dimensional invertible representation $(\mathbf{V}',\bs\phi')$ satisfying $\mu'(\bs\phi')=\bs\zeta'$, where $\bs\zeta'$, $\mathbf{d}'$ and $\bs\theta'$ are defined by reflection, and $\mathbf{V}'$ and $\bs\phi'$ are defined as follows. 

The definition of $(\mathbf{V}',\bs\phi')$ is local at $v$. Without loss of generality, we may assume that there is no arrow $a\in Q_1$ with $t(a)=v$. We then define $Q_1(v):=\{a\in Q_1\mid h(a)=v\}$, $Q_0(v):=\{v\}\cup\{t(a)\mid a\in Q_1(v)\}$, $Q(v)=(Q_0(v),Q_1(v))$ and $\bar{Q}(v):=\overline{Q(v)}$. We also write $\hat{V}_v:=\oplus_{a\in Q_1(v)}V_{t(a)}$. For any $a\in Q_1(v)$, let $\pi_a:\hat{V}_v\rightarrow V_{t(a)}$ and $\iota_a:V_{t(a)}\rightarrow\hat{V}_v$ be the natural projection and inclusion respectively, and define $\Phi_a:V_v\rightarrow V_v$ by
$$
\Phi_a:=\prod^{\prec}_{\substack{a'\in Q_1(v)\\a'\prec a}}(1+\phi_{a'}\phi_{a^{\prime\ast}}).
$$
We then define
$$
\iota_v=\Sigma_{a\in Q_1(v)}\iota_a\phi_{a^{\ast}}:V_v\rightarrow\hat{V}_v\text{ and }\pi_v=\Sigma_{a\in Q_1(v)}\Phi_a\phi_a\pi_a:\hat{V}_v\rightarrow V_v. 
$$
Beware that the map $\pi_v$ differs from the one defined by Crawley-Boevey-Shaw by a scalar so that it makes sense when $\zeta_v=1$. Then, the definition of $(\mathbf{V}',\bs\phi')$ is
$$
V'_v=\Ker\pi_v \text{ and } V'_w=V_w\text{ for }w\ne v,
$$
with $\bs\phi'=(\phi'_a)_{a\in Q_1\sqcup Q_1^{\ast}}$ defined by
\begingroup
\allowdisplaybreaks
\begin{align*}
\phi'_a&=\phi_a\text{ if $a$ is not incident at $v$},\\
\phi'_a&=\pi_{a^{\ast}}|_{V'_v}:V'_v\rightarrow V_{h(a)}\text{ if }t(a)=v,\text{ and}\\
\phi_a'&=\sum_{\substack{a'\in Q_1(v)\\a'\prec a}}\iota_{a'}\phi_{a^{\prime\ast}}\phi_a+\frac{1}{q_v}\sum_{\substack{a'\in Q_1(v)\\a'\succeq a}}\iota_{a'}\phi_{a^{\prime\ast}}\phi_a+\frac{1-q_v}{q_v}\iota_a:V_{t(a)}\rightarrow V'_v\text{ if }h(a)=v.
\end{align*}
\endgroup

For an arrow $a$ with $h(a)=v$, we only know a priori that $\phi_a'$ maps to $\hat{V}_v$. The proofs of \cite[Lemma 4.1]{CBS} and \cite[Lemma 4.2]{CBS} remain valid in the case $\zeta_v=1$, as long as we put the factor $\zeta_v-1$ in the right place. We thus see that $\phi_a'$ has its image in $V'_v$. It also follows that $\mu'_v(\bs\phi')=\zeta'_v$ even if $\zeta_v=1$. It was also explained in \cite[\S 5.1]{Yama} that $\mu'_w(\bs\phi')=\zeta'_w$ for $w\ne v$. The claim that $\mathbf{V}'$ has dimension $\mathbf{d}'$ follows from the surjectivity of $\pi_v$; while the case $\zeta_v\ne 1$ was explained by Yamakawa, the case $\zeta_v=1$ follows from the lemma below.
\begin{Lem}\label{Lem-Yama-3.7}
Let $\bs\phi\in\Rep^{\circ}(\bar{Q},\mathbf{d})^{\bs\theta-ss}$. Then, we have:
\begin{itemize}
\item[(i)] If $\theta_v>0$, then $\iota_v$ is injective.
\item[(ii)] If $\theta_v<0$, then $\pi_v$ is surjective.
\end{itemize}
\end{Lem}
\begin{proof}
This is a mild modification of \cite[Proposition 3.7]{Yama}.
\end{proof}
If $\theta_v<0$, then the representations $(\mathbf{V},\bs\phi)$ and $(\mathbf{V}',\bs\phi')$ satisfy the assumptions of \cite[Lemma 5.4]{Yama}, and the same arguments for this lemma show that $\bs\phi$ is $\bs\theta$-semi-stable if and only if $\bs\phi'$ is $\bs\theta'$-semi-stable. Therefore, $(\mathbf{V}',\bs\phi')$ indeed defines an element of $\Rep^{\circ}(\bar{Q},\bs\zeta',\mathbf{d}')^{\bs\theta'-ss}$. The rest of this subsection will be devoted to showing that this construction induces an isomorphism between multiplicative quiver varieties.

As above, we consider the subquiver $Q(v)$ associated to $v$.  Let $\mathbf{d}_1$ and $\mathbf{d}_2$ be dimension vectors in $\mathbb{Z}^{Q_0(v)}$ satisfying
\begin{itemize}
\item $d_{1,w}=d_{2,w}=d_w$ for any $w\ne v$ and some $d_w\in\mathbb{Z}_{>0}$, and
\item $d_{1,v}+d_{2,v}=\sum_{w\ne v}d_w$.
\end{itemize} 
We will write $\hat{d}_v=d_{1,v}+d_{2,v}$ in what follows.
\begin{Lem}\label{Lem-S(v)}
Let $\zeta_v\in\mathbb{C}^{\ast}$. Consider the locally closed subvariety
$$
\mathscr{S}(v)\subset \Rep(\bar{Q}(v),\mathbf{d}_1)\times\Rep(\bar{Q}(v),\mathbf{d}_2)
$$
consisting of the pairs of quiver representations $(\bs\phi,\bs\psi)$ satisfying the following conditions:
\begin{itemize}
\item[(1a)] $\sum_{a\in Q_1(v)}\phi_a\phi_{a^{\ast}}=\zeta_v-1$,
\item[(1b)] $\sum_{a\in Q_1(v)}\psi_a\psi_{a^{\ast}}=\zeta_v^{-1}-1$,
\item[(2)] $(\Sigma_{a\in Q_1(v)}\iota_a\phi_{a^{\ast}})(\Sigma_{a\in Q_1(v)}\phi_a\pi_a)=\zeta_v(\Sigma_{a\in Q_1(v)}\iota_a\psi_{a^{\ast}})(\Sigma_{a\in Q_1(v)}\psi_a\pi_a)+\zeta_v-1$,
\item[(3a)] $\bar{\phi}:=\Sigma_{a\in Q_1(v)}\phi_a\pi_a$ is surjective,
\item[(3b)] $\bar{\psi}^{\ast}:=\Sigma_{a\in Q_1(v)}\iota_a\psi_{a^{\ast}}$ is injective; and
\item[(4)] the sequence
$$
\mathbb{C}^{d_{2,v}}\stackrel{\bar{\psi}^{\ast}}{\longrightarrow} \mathbb{C}^{\hat{d}_v}\stackrel{\bar{\phi}}{\longrightarrow}\mathbb{C}^{d_{1,v}}
$$
is exact.
\end{itemize}
Define the locally closed subvariety
$$
\mathscr{X}(v):=\{\bs\phi\in\Rep(\bar{Q}(v),\mathbf{d}_1)\mid \text{(1a) and (3a) hold}\}.
$$
Then, the fibre of the projection $\mathscr{S}(v)\rightarrow\mathscr{X}(v)$ is nonempty and is isomorphic to $\GL_{d_{2,v}}$.
\end{Lem}
\begin{proof}
We will write $\bar{\phi}^{\ast}:=\Sigma_{a\in Q_1(v)}\iota_a\phi_{a^{\ast}}$, which is a map $\mathbb{C}^{d_{1,v}}\rightarrow\mathbb{C}^{\hat{d}_v}$. With this notation, condition (1a) says $\bar\phi\bar\phi^{\ast}=\zeta_v-1$ as an endomorphism of $\mathbb{C}^{d_{1,v}}$. The quiver representation $\bs\phi$ is equivalent to a pair of linear maps $(\bar\phi,\bar\phi^{\ast})$, and similarly for $\bs\psi$. Then, condition (2) says that $\bar\phi^{\ast}\bar\phi=\zeta_v\bar\psi^{\ast}\bar\psi+\zeta_v-1$.

Suppose that $\zeta_v\ne 1$. Note that (3a) and (3b) are consequences of (1a) and (1b) in this case. It follows from the proof of \cite[Lemma 2.1]{CB01} that the fibre over $\bs\phi\in\mathscr{X}(v)$ is isomorphic to $\GL_{d_{2,v}}$. Indeed, the condition (4) also follows from (1a), (1b) and (2), and thus it is redundant.

Now, suppose $\zeta_v=1$ and let us compute the fibre over a given $\bs\phi\in\mathscr{X}(v)$. Let $K:=\Ker\bar\phi$, which is a $d_{2,v}$-dimensional space. Let $\bar\psi^{\ast}$ be any isomorphism from $\mathbb{C}^{d_{2,v}}$ to $K$. In view of (1a), we may regard $\bar\phi^{\ast}\bar\phi$ as a map $\mathbb{C}^{\hat{d}_v}\rightarrow K$. Then, (2) implies that $\bar\psi$ is completely determined by $\bar\psi^{\ast}$ and $\bar\phi^{\ast}\bar\phi$. Now, we check that such a pair $(\bar\psi,\bar\psi^{\ast})$ also satisfies (1b). Since $\bar\psi^{\ast}$ is injective, it suffices to show that $\bar\psi^{\ast}\bar\psi\bar\psi^{\ast}=0$, but $ \bar\psi^{\ast}\bar\psi\bar\psi^{\ast}=\bar\phi^{\ast}\bar\phi\bar\psi^{\ast}$ and $\bar\phi\bar\psi^{\ast}=0$. Therefore, the fibre is nonempty and is isomorphic to $\GL_{d_{2,v}}$.
\end{proof}

\begin{Prop}(Yamakawa)\label{Prop-ref-isom}
If $v$ is a loop-free vertex such that either $\theta_v\ne 0$ or $\zeta_v\ne 1$, then the reflection functor induces an isomorphism 
\begin{equation}\label{eq-Prop-ref-isom}
\mathcal{M}(\bs\zeta,\mathbf{d},\bs\theta)\lisom\mathcal{M}(u_v(\bs\zeta),s_v(\mathbf{d}),r_v(\bs\theta)).
\end{equation}
\end{Prop}
\begin{proof}
We need to check that the arguments for \cite[Theorem 5.3]{Yama} remain valid when the strictly semi-stable loci are nonempty. Write $\bs\zeta'=u_v(\bs\zeta)$, $\mathbf{d}'=s_v(\mathbf{d})$ and $\bs\theta'=r_v(\bs\theta)$. Define the correspondence 
$$
Z(\bs\zeta,\mathbf{d},\bs\theta)\subset\Rep^{\circ}(\bs\zeta,\mathbf{d})^{\bs\theta-ss}\times\Rep^{\circ}(\bs\zeta',\mathbf{d}')^{\bs\theta'-ss}
$$
to be the closed subvariety consisting of pairs $((\mathbf{V},\bs\phi),(\mathbf{V}'
,\bs\phi'))$ satisfying
\begin{itemize}
\item[(1)] $\phi_a=\phi_a'$ for any arrow $a$ not incident at $v$,
\item[(2)] the sequence 
$$
0\longrightarrow V_v'\stackrel{\iota_v'}{\longrightarrow} \hat{V}_v\stackrel{\pi_v}{\longrightarrow} V_v\longrightarrow 0
$$
is exact, where $\iota_v'$ is defined by $\bs\phi'$; and
\item[(3)] $\iota_v\pi_v=\zeta_v\iota_v'\pi_v'+\zeta_v-1$.
\end{itemize}
Note that the injectivity and surjectivity in (2) are open conditions. If $\zeta_v\ne1$, then the relation $\mu_v(\bs\phi)=1+\pi_v\iota_v=\zeta_v$ implies that $\pi_v$ is surjective, and the injectivity of $\iota'_v$ follows similarly. If $\zeta_v=1$, then we may assume that $\theta_v<0$ and $\theta'_v>0$, in which case Lemma \ref{Lem-Yama-3.7} implies the injectivity and surjectivity in (2). The open conditions are thus vacuous and $Z(\bs\zeta,\mathbf{d},\bs\theta)$ is indeed a closed subvariety. (We will not need the closedness in the proof though.)

Let $p_1:Z(\bs\zeta,\mathbf{d},\bs\theta)\rightarrow\Rep^{\circ}(\bs\zeta,\mathbf{d})^{\bs\theta-ss}$ (resp. $p_2:Z(\bs\zeta,\mathbf{d},\bs\theta)\rightarrow\Rep^{\circ}(\bs\zeta',\mathbf{d}')^{\bs\theta'-ss}$) be the restriction of the projection to the first (resp. second) factor. We will show that $p_1$ is a principal $\GL_{d'_v}$-bundle and that $p_2$ is a principal $\GL_{d_v}$-bundle. It follows that they induce isomorphisms
$$
\Rep^{\circ}(\bs\zeta,\mathbf{d})^{\bs\theta-ss}\ds\GL_{\mathbf{d}}\lisomL Z(\bs\zeta,\mathbf{d},\bs\theta)\ds\GL_{\tilde{\mathbf{d}}}\lisom\Rep^{\circ}(\bs\zeta',\mathbf{d}')^{\bs\theta'-ss}\ds\GL_{\mathbf{d}'},
$$
as desired. We will only give the details for the statement that $p_1$ is surjective with fibres isomorphic to $\GL_{d_v'}$, since Yamakawa has shown that it admits local sections in \cite[\S 5.4]{Yama}, and the arguments are valid as long as we have the surjectivity of $\pi_v$.

Consider the quiver $Q(v)$. As above, we assume that there is no arrow $a\in Q_1(v)$ with $t(a)=v$. Let $\mathbf{d}_1$ (resp. $\mathbf{d}_2$) be the restriction of $\mathbf{d}$ (resp. $\mathbf{d}'$) to $Q(v)$. Let 
$$
\mathscr{R}(v)\subset\Rep^{\circ}(\bar{Q}(v),\mathbf{d}_1)\times\Rep^{\circ}(\bar{Q}(v),\mathbf{d}_2)
$$
be the subvariety consisting of the pairs of quiver representations $(\bs\phi,\bs\phi')$ satisfying 
\begin{itemize}
\item[(1a)] $\mu_v(\bs\phi)=\zeta_v$,
\item[(1b)] $\mu'_v(\bs\phi')=\zeta_v^{-1}$,
\item[(2)] $\iota_v\pi_v=\zeta_v\iota_v'\pi_v'+\zeta_v-1$; and
\item[(4)] the sequence
$$
0\longrightarrow\mathbb{C}^{d_{2,v}}\stackrel{\iota_v'}{\longrightarrow} \mathbb{C}^{\hat{d}_v}\stackrel{\pi_v}{\longrightarrow}\mathbb{C}^{d_{1,v}}\longrightarrow0
$$
is exact.
\end{itemize}
Let $\Rep^{\circ}(\bar{Q}(v),\zeta_v,\mathbf{d}_1)\subset\Rep^{\circ}(\bar{Q}(v),\mathbf{d}_1)$ be the subvariety defined by (1a) and the condition that $\pi_v$ be surjective. 
For $i\in\{1,2\}$, let $U_i\subset\Rep(\bar{Q}(v),\mathbf{d}_i)$ be the open subset consisting of quiver representations $\tilde{\bs\phi}$ satisfying
$$
\det(1+\sum_{\substack{a'\in Q_1(v)\\a'\preceq a}}\tilde{\phi}_{a'}\tilde{\phi}_{a^{\prime\ast}})\ne0
$$
for any $a\in Q_1(v)$. Then, there is an isomorphism
\begingroup
\allowdisplaybreaks
\begin{align*}
f_i:\Rep^{\circ}(\bar{Q}(v),\mathbf{d}_i)&\longrightarrow U_i\\
\bs\phi&\longmapsto \tilde{\bs\phi}=(\tilde{\phi}_a,\tilde{\phi}_{a^{\ast}})_{a\in Q_1(v)}
\end{align*}
\endgroup
with $\tilde{\phi}_{a^{\ast}}=\phi_{a^{\ast}}$ and 
$$
\tilde{\phi}_a=\prod^{\prec}_{\substack{a'\in Q_1(v)\\a'\prec a}}(1+\phi_{a'}\phi_{a^{\prime\ast}})\cdot\phi_a.
$$

In Lemma \ref{Lem-S(v)}, we have introduced the varieties $\mathscr{S}(v)$ and $\mathscr{X}(v)$. Let $\mathscr{S}^{\circ}(v):=\mathscr{S}(v)\cap U_1\times U_2$; then, $f_1\times f_2$ restricts to an isomorphism $f$ which fits into the following commutative diagram
\begin{equation}\label{eq-red-diag}
\begin{tikzcd}[row sep=2.5em, column sep=2.5em]
\mathscr{R}(v)\arrow[r, "f"] \arrow[d, "\pr_1"'] & \mathscr{S}^{\circ}(v) \arrow[d, "\pr'_1"]\\
\Rep^{\circ}(\bar{Q}(v),\mathbf{d}_1) \arrow[r, "f_1"'] & U_1.
\end{tikzcd}
\end{equation}
Let $\mathscr{X}^{\circ}(v):=\mathscr{X}(v)\cap U_1$. We claim that $\pr'_1$ surjects onto $\mathscr{X}^{\circ}(v)$ and has fibres isomorphic to $\GL_{d_{2,v}}$. For any $\tilde{\bs\phi}\in\mathscr{X}^{\circ}(v)$, we know that its preimage in $\mathscr{S}(v)$ is isomorphic to $\GL_{d_{2,v}}$; we need to show that this preimage is contained in $\mathscr{S}^{\circ}(v)$. Let $\bs\phi_{(v)}\in\Rep^{\circ}(\bar{Q}(v),\mathbf{d}_1)$ be such that $f_1(\bs\phi_{(v)})=\tilde{\bs\phi}$. Let $\bs\phi_{(v)}'\in\Rep^{\circ}(\bar{Q}(v),\mathbf{d}_2)$ be determined by $\bs\phi_{(v)}$ as in the beginning of this subsection. If $\zeta_v\ne1$, then the proof of \cite[Lemma 4.3]{CBS} shows that $\bs\phi_{(v)}'$ is invertible at $v$. If $\zeta_v=1$, then similar computations also show that $\bs\phi_{(v)}'$ is invertible at $v$. Now, we have $(\bs\phi_{(v)},\bs\phi_{(v)}')\in\mathscr{R}(v)$ by \cite[Proposition 5.6]{Yama}, and so its image under $f$ lies over $\tilde{\bs\phi}$; it follows that the fibre of $\pr_1'$ over $\tilde{\bs\phi}$ is nonempty. It is easy to see from the definition that $U_2$ is invariant under the action of $\GL_{d_{2,v}}$. It follows that $\mathscr{S}^{\circ}(v)\rightarrow\mathscr{X}^{\circ}(v)$ is surjective and has fibres isomorphic to $\GL_{d_{2,v}}$, and the same is true for the map $\mathscr{R}(v)\rightarrow\Rep^{\circ}(\bar{Q}(v),\zeta_v,\mathbf{d}_1)$.

Let us show that
$$
p_1:Z(\bs\zeta,\mathbf{d},\bs\theta)\rightarrow\Rep^{\circ}(\bs\zeta,\mathbf{d})^{\bs\theta-ss}
$$
is surjective with fibres isomorphic to $\GL_{d_{2,v}}$. Let $Q[v]$ be the quiver obtained from $Q$ by removing the vertex $v$ and the arrows incident at $v$. Then, $\Rep^{\circ}(\bs\zeta,\mathbf{d})^{\bs\theta-ss}$ is the subvariety of
$$
\Rep^{\circ}(\bar{Q}(v),\zeta_v,\mathbf{d}_1)\times \Rep^{\circ}(\bar{Q}[v],\mathbf{d}|_{Q[v]})\subset\Rep^{\circ}(\bar{Q},\mathbf{d})
$$ 
defined by the equations $\mu_w(\bs\phi)=\zeta_w$ for $w\ne v$ and the stability condition; here, the surjectivity assumption in the definition of $\Rep^{\circ}(\bar{Q}(v),\zeta_v,\mathbf{d}_1)\subset\Rep^{\circ}(\bar{Q}(v),\mathbf{d}_1)$ is always satisfied since either $\zeta_v\ne 1$ or $\theta_v<0$. However, $Z(\bs\zeta,\mathbf{d},\bs\theta)$ is the subvariety of
$$
\mathscr{R}(v)\times\Rep^{\circ}(\bar{Q}[v],\mathbf{d}|_{Q[v]})
$$
defined by the stability condition and the equations $\mu_w(\bs\phi)=\zeta_w$ for $w\ne v$, where for any $(\bs\phi_1,\bs\phi_2)\in\mathscr{R}(v)$ and $\bs\phi_3\in\Rep^{\circ}(\bar{Q}[v],\mathbf{d}|_{Q[v]})$ we write $\bs\phi=(\bs\phi_1,\bs\phi_3)\in\Rep^{\circ}(\bar{Q},\mathbf{d}_1)$. Note that $\mu_w(\bs\phi)=\zeta_w$ is equivalent to $\mu'_w(\bs\phi')=\zeta'_w$ for $w\ne v$ if we write $\bs\phi'=(\bs\phi_2,\bs\phi_3)$; this follows from the condition (2) in the definition of $\mathscr{R}(v)$. Besides, it follows from the proof of \cite[Lemma 5.4]{Yama} that $\bs\phi$ is $\bs\theta$-semi-stable if and only if $\bs\phi'$ is $\bs\theta'$-semi-stable. Therefore, $Z(\bs\zeta,\mathbf{d},\bs\theta)$ is the inverse image of $\Rep^{\circ}(\bs\zeta,\mathbf{d})^{\bs\theta-ss}$ under the following map
$$
\mathscr{R}(v)\times\Rep^{\circ}(\bar{Q}[v],\mathbf{d}|_{Q[v]})\longrightarrow\Rep^{\circ}(\bar{Q}(v),\zeta_v,\mathbf{d}_1)\times \Rep^{\circ}(\bar{Q}[v],\mathbf{d}|_{Q[v]}),
$$ 
which is surjective with fibres isomorphic to $\GL_{d_{2,v}}$. 
\end{proof}

\subsection{Multiplicative quiver varieties for type A quivers}\label{subsec-multi-Q-typeA}\hfill 

For now, we write $G_0=\GL_n$. Let $P_0\subset G_0$ be a proper parabolic subgroup. The parabolic subgroups conjugate to $P_0$ are parametrised by $G_0/P_0$. Let $L_0$ be the Levi quotient of $P_0$ and $U_0$ the unipotent radical of $P_0$. Let $\bar{C}$ be the closure of an $L_0$-conjugacy class $C\subset L_0$. We introduce the following variety
\begin{equation}\label{eq-defn-F(C)}
\mathbf{F}(\bar{C}):=\{(X,gP_0)\in G\times G/P_0\mid g^{-1}Xg\in\bar{C}U_0\},
\end{equation}
where $\bar{C}U_0$ is the inverse image of $\bar{C}$ in $P_0$. There is an action of $G_0$ on $\mathbf{F}(\bar{C})$ given by $h:(X,gP_0)\mapsto(hXh^{-1},hgP_0)$ for any $h\in G_0$. The projection of $\mathbf{F}(\bar{C})$ to $G_0$ is the closure of a conjugacy class $C^{\sharp}$ (with its closure denoted by $\bar{C}^{\sharp}$). We thus see that $\mathbf{F}(\bar{C})$ is projective over an affine variety.
\begin{Rem}\label{Rem-CaEs-paper}
These assertions are equivalently stated as \cite[Lemma 4.3]{CaEs}. The key fact is that $G_0.(CU_0)$, the union of the $G_0$-conjugates of $CU_0$, meets only finitely many conjugacy classes in $G_0$.
\end{Rem}

Our goal is to realise $\mathbf{F}(\bar{C})$ as a multiplicative quiver variety. We define a type $A$ quiver from $\mathbf{F}(\bar{C})$ as follows. We first regard $P_0$ as the stabiliser of a flag $F_{\bullet}$:
\begin{equation}\label{eq-subsec-Ch-var-multi-Q-1}
\mathbb{C}^n=F_0\supset F_1\supset\cdots \supset F_{\mu}\supset F_{\mu+1}=0.
\end{equation}
Write $\Gr F_{\bullet}:=\bigoplus_{s=0}^{\mu}F_s/F_{s+1}$ and $\GL(\Gr F_{\bullet}):=\prod_{s=0}^{\mu}\GL(F_s/F_{s+1})$. We then regard $\bar{C}$ as the closure of a conjugacy class of $\GL(\Gr F_{\bullet})$. For each $s$, we denote by $\bar{C}^{(s)}\subset\GL(F_s/F_{s+1})$ the corresponding component of $\bar{C}$, and let $w_s+1$ be the degree of the minimal polynomial of an element of $C^{(s)}$. For $0\le s\le \mu$, we write $i_s=\sum_{0\le s'<s}(w_{s'}+1)$ and write $i_{\mu+1}=\nu+1=\sum_{0\le s\le\mu}(w_s+1)$. We now define the double quiver $\bar{Q}=(Q_0,\bar{Q}_1)$, where $Q_0=\{0,1,\ldots,\nu\}$ and 
$$
\bar{Q}_1=\{a_i:i\mapsto i+1,a_i^{\ast}:i+1\mapsto i\mid 0\le i\le \nu-1\}.$$
The vertex 0 will be regarded as a framing vertex.

Next, we continue to define a quiver variety from $\mathbf{F}(\bar{C})$. Let $(X,gP_0)\in\mathbf{F}(\bar{C})$. Since $g^{-1}Xg$ preserves the flag $F_{\bullet}$, it defines an element $X^{(s)}\in\GL(F_s/F_{s+1})$ for every $0\le s\le\mu$. By the definition of $\mathbf{F}(\bar{C})$, we have $X^{(s)}\in\bar{C}^{(s)}$. We may choose $(X,gP_0)$ so that $X^{(s)}\in C^{(s)}$ for every $s$. Then, we can apply the construction in the trivial stability case to obtain a flag on $F_s/F_{s+1}$ as follows; see, e.g., \cite[\S 1]{CB04}. Let $(\xi^{(s)}_t)_{0\le t\le w_s}$ be a tuple of complex numbers satisfying 
$$
\prod_{0\le t\le w_s}(M^{(s)}-\xi^{(s)}_t\Id)=0.
$$
For $1\le t\le w_s$, define
$$
W^{(s)}_t:=\Ima (X^{(s)}-\xi^{(s)}_0)(X^{(s)}-\xi^{(s)}_1)\cdots(X^{(s)}-\xi^{(s)}_{t-1})\subset F_s/F_{s+1}.
$$
These subspaces form a flag:
$$
F_s/F_{s+1}=W^{(s)}_0\supset W^{(s)}_1\supset \cdots\supset W^{(s)}_{w_s}\supset W^{(s)}_{w_s+1}=0.
$$
For $0\le s\le \mu$, $0\le t\le w_s$ and $i=i_s+t$, we denote by $F_i\subset V_0$ the subspace containing $F_{s+1}$ satisfying $F_i/F_{s+1}=W^{(s)}_t$. We thus obtain a flag $V_{\bullet}$:
\begin{equation}\label{eq-defn-V}
F_0=V_0\supset V_1\supset\cdots\supset V_{\nu}\supset V_{\nu+1}=0,
\end{equation}
which refines $F_{\bullet}$. For $0\le s\le\mu$, we have $V_{i_s}=F_s$. Let $\mathbf{d}=(d_i)_{i\in Q_0}$ with $d_i=\dim V_i$. Define $\bs\theta=(\theta_i)_{i\in Q_0}$ in such a way that $\theta_i\in\mathbb{Q}_{>0}$ whenever $i=i_s$ for some $1\le s\le \mu$ and $\theta_i=0$ otherwise. Define $\bs\xi=(\xi_i)_{i\in Q_0}$ in such a way that for any $0\le s\le\mu$ and $i_s\le i<i_{s+1}$, we have $\xi_i=\xi^{(s)}_{i-i_s}$. We then define $\bs\zeta=(\zeta_i)_{i\in Q_0}$ so that $\zeta_0=\xi_0$ and $\zeta_i\xi_{i-1}=\xi_i$ for any $i>0$. By abuse of notation, we will also write $\bs\zeta=(\zeta_i)_{i\in Q_0\setminus\{0\}}$ if only the vertices $i>0$ are relevant.

Define the map
\begingroup
\allowdisplaybreaks
\begin{align*}
\mu_{>0}:\Rep^{\circ}(\bar{Q},\mathbf{d})&\longrightarrow\prod_{i\in Q_0\setminus\{0\}}\GL_{d_i}\\
\bs\phi&\longmapsto ((\Id+\phi_{i-1}\phi^{\ast}_{i-1})(\Id+\phi^{\ast}_i\phi_i)^{-1})_{i\in Q_0\setminus\{0\}},
\end{align*}
\endgroup
where $\phi_i:=\phi_{a_i}$ and $\phi_i^{\ast}:=\phi_{a_i^{\ast}}$. The fibre $\mu_{>0}^{-1}(\bs\zeta)$ is a $G_0$-invariant subvariety, with the group $G_0$ naturally acting on $F_0$. The relevant stability condition, however, is not the usual one defined for the representations of multiplicative preprojective algebras. It is instead the stability condition for framed quiver representations. Quiver varieties for framed quivers were first studied by Nakajima in \cite{Naka94} and \cite{Naka98}, and some properties are nicely summarised in \cite{Naka09}. Let $(\mathbf{V},\bs\phi)$ be a representation of $\bar{Q}$ with dimension vector $\mathbf{d}$. We say that $(\mathbf{V},\bs\phi)$ is $\bs\theta$-stable (resp. $\bs\theta$-semi-stable) if for any subrepresentation $(\mathbf{V}',\bs\phi'=\bs\phi|_{\mathbf{V}'})$ with $V'_0=\{0\}$ and $\dim V'_i=d'_i$, we have
\begin{equation}\label{eq-semi-st-1}
\sum_{i\in Q_0\setminus\{0\}}d'_i\theta_i<0\text{ (resp. $\le0$)},
\end{equation}
and for any subrepresentation $(\mathbf{V}',\bs\phi'=\bs\phi|_{\mathbf{V}'})$ with $V'_0=V_0$ and $\dim V'_i=d'_i$, we have
\begin{equation}\label{eq-semi-st-2}
\sum_{i\in Q_0\setminus\{0\}}d'_i\theta_i<\sum_{i\in Q_0\setminus\{0\}}d_i\theta_i\text{ (resp. $\le\sum_{i\in Q_0\setminus\{0\}}d_i\theta_i$)}
\end{equation}
(beware that the vertex $0$ is a framing vertex).
The open subset of $\bs\theta$-semi-stable (resp. $\bs\theta$-stable) representations in $\Rep^{\circ}(\bar{Q},\mathbf{d})$ will be denoted by $\Rep^{\circ}(\bar{Q},\mathbf{d})^{\bs\theta-ss}$ (resp. $\Rep^{\circ}(\bar{Q},\mathbf{d})^{\bs\theta-s}$). We prepare some lemmas before proving Proposition \ref{Prop-F(C)-unframed-quiver}.

\begin{Lem}\label{Lem-Ker-subrepres}
Let $(\mathbf{V},\bs\phi)$ be a quiver representation representing a point of $\mu_{>0}^{-1}(\bs\zeta)$. Let $m\in Q_0\setminus\{0\}$. Define a $Q_0$-graded vector space $\mathbf{W}$ by 
$$
W_{m}=\Ker\phi_{0}^{\ast}\cdots\phi_{m-1}^{\ast},
$$ 
$W_{i}:=\phi_{i}^{\ast}\cdots\phi_{m-1}^{\ast}(W_{m})$ if $i<m$, and $W_{i}:=\phi_{i-1}\cdots\phi_{m}(W_{m})$ if $i>m$; in particular, we have $W_0=0$. Then, $\mathbf{W}$ is a subrepresentation of $(\mathbf{V},\bs\phi)$.
\end{Lem}
\begin{proof}
We need to show that $\mathbf{W}$ is preserved by $\bs\phi$. As the first step, we show that $\phi_{m-1}(W_{m-1})\subset W_{m}$, or equivalently, that for any $v\in W_{m}$, we have
$$
\phi_{0}^{\ast}\cdots\phi_{m-1}^{\ast}\phi_{m-1}\phi_{m-1}^{\ast}(v)=0.
$$
We compute
\begingroup
\allowdisplaybreaks
\begin{align*}
\phi_{0}^{\ast}\cdots\phi_{m-1}^{\ast}\phi_{m-1}\phi_{m-1}^{\ast}(v)&=\phi_{0}^{\ast}\cdots\phi_{m-2}^{\ast}(1+\phi_{m-1}^{\ast}\phi_{m-1})\phi_{m-1}^{\ast}(v)\\
&=\zeta_{m-1}^{-1}\phi_{0}^{\ast}\cdots\phi_{m-2}^{\ast}(1+\phi_{m-2}\phi_{m-2}^{\ast})\phi_{m-1}^{\ast}(v)\\
&=\zeta_{m-1}^{-1}\phi_{0}^{\ast}\cdots\phi_{m-2}^{\ast}\phi_{m-2}\phi_{m-2}^{\ast}\phi_{m-1}^{\ast}(v).
\end{align*}
\endgroup
By induction, we have
$$
\phi_{0}^{\ast}\cdots\phi_{m-1}^{\ast}\phi_{m-1}\phi_{m-1}^{\ast}(v)=\big(\prod_{i=1}^{m-1}\zeta_{i}^{-1}\big)\phi_{0}^{\ast}\phi_{0}\phi_{0}^{\ast}\cdots\phi_{m-1}^{\ast}(v)=0.
$$
We have thus proved the inclusion $\phi_{m-1}(W_{m-1})\subset W_{m}$.

Next, we show that 
\begin{equation}\label{eq-Lem-ss->inj}
\phi_{i-1}(W_{i-1})\subset W_{i}
\end{equation}
for any $i<m$. Suppose that the inclusion has been proved for $\phi_{i}$. Let us show that for any $v\in W_{i}$, we have $\phi_{i-1}\phi_{i-1}^{\ast}(v)\in W_{i}$. We have,
$$
\phi_{i-1}\phi_{i-1}^{\ast}(v)=(1+\phi_{i-1}\phi_{i-1}^{\ast})(v)-v=\zeta_{i}(1+\phi_{i}^{\ast}\phi_{i})(v)-v.
$$
However, $\phi_{i}^{\ast}\phi_{i}(v)\in\phi_{i}^{\ast}(W_{i+1})=W_{i}$ by the induction hypothesis. It follows that the desired inclusion holds for all $i<m$.

We then show that $\phi_{i-1}^{\ast}(W_{i})\subset W_{i-1}$ for any $i>m$. This will follow from the following equation
$$
\xi_{i-1}(1+\phi^{\ast}_{i-1}\phi_{i-1})\phi_{i-2}\cdots\phi_{m}=\xi_{m-1}\phi_{i-2}\cdots\phi_{m}(1+\phi_{m-1}\phi^{\ast}_{m-1})
$$
which appeared in the proof of \cite[Lemma 4.3]{Yama}. For the convenience of the reader, let us explain how it is proved. Recall that $\zeta_{k}=\xi_{k}/\xi_{k-1}$ for any $k>0$. We compute
\begingroup
\allowdisplaybreaks
\begin{align*}
\xi_{m-1}\phi_{i-2}\cdots\phi_{m}(1+\phi_{m-1}\phi^{\ast}_{m-1})&=\zeta_{m}\xi_{m-1}\phi_{i-2}\cdots\phi_{m}(1+\phi^{\ast}_{m}\phi_{m})\\
&=\xi_{m}\phi_{i-2}\cdots\phi_{m+1}(1+\phi_{m}\phi^{\ast}_{m})\phi_{m}.
\end{align*}
\endgroup
By inductin, we have
\begingroup
\allowdisplaybreaks
\begin{align*}
\xi_{m-1}\phi_{i-2}\cdots\phi_{m}(1+\phi_{m-1}\phi^{\ast}_{m-1})&=\big(\prod_{i'=m}^{i-2}\zeta_{i'}\big)\xi_{m-1}(1+\phi_{i-2}\phi^{\ast}_{i-2})\phi_{i-2}\cdots\phi_{m}\\
&=\xi_{i-1}(1+\phi^{\ast}_{i-1}\phi_{i-1})\phi_{i-2}\cdots\phi_{m}.
\end{align*}
\endgroup
This completes proving that $\mathbf{W}$ is a subrepresentation.
\end{proof}

\begin{Lem}\label{Lem-ss->inj}
Let $(\mathbf{V},\bs\phi)$ be a $\bs\theta$-semi-stable representation representing a point of 
$$
\Rep^{\circ}(\bar{Q},\mathbf{d})^{\bs\theta-ss}\cap\mu_{>0}^{-1}(\bs\zeta),
$$ 
and let $m=i_s$ for some $1\le s\le \mu$ so that $\theta_m>0$. Then, the following assertions hold:
\begin{itemize}
\item[(i)] The map $\phi_{0}^{\ast}\cdots\phi_{m-1}^{\ast}$ is injective.
\item[(ii)] The subspace $F_s:=\Ima\phi_{0}^{\ast}\cdots\phi_{m-1}^{\ast}$ is preserved by $\zeta_0(\Id+\phi^{\ast}_{0}\phi_{0})$.
\end{itemize}
\end{Lem}

\begin{proof}
Consider the subrepresentation $\mathbf{W}\subset \mathbf{V}$ associated to $m$ as in the above lemma. If the asserted injectivity of (i) fails, then this subrepresentation destabilises $V$, contradicting our assumption. In order to show that $\zeta_{0}(1+\phi^{\ast}_{0}\phi_{0})$ preserves $F_s$, we consider the subspaces $W'_{m}:=V_{m}$ and $W'_{i}:=\phi_{i}^{\ast}\cdots\phi_{m-1}^{\ast}(W'_{i})$ for $i<m$. Then, equation (\ref{eq-Lem-ss->inj}) holds for $W'_i$ for all $i\le m$, which follows from the same computation. Therefore, $\zeta_{0}(1+\phi^{\ast}_{0}\phi_{0})$ preserves $W'_0=F_s$.
\end{proof}

\begin{Lem}\label{Lem-st->inj}
Let $(\mathbf{V},\bs\phi)$ be a $\bs\theta$-stable representation representing a point of 
$$
\Rep^{\circ}(\bar{Q},\mathbf{d})^{\bs\theta-s}\cap\mu_{>0}^{-1}(\bs\zeta). 
$$
Let $m\in Q_0$ be such that $i_s<m<i_{s+1}$ for some $0\le s\le \mu$. Then, the map $\phi_{i_s}^{\ast}\cdots\phi_{m-1}^{\ast}$ is injective.
\end{Lem}
\begin{proof}
Consider the subrepresentation $\mathbf{W}\subset \mathbf{V}$ associated to $m$. Since $\phi_0^{\ast}\cdots\phi_{i_s-1}^{\ast}$ is injective by Lemma \ref{Lem-ss->inj}, the kernel of $\phi_{i_s}^{\ast}\cdots\phi_{m-1}^{\ast}$ is precisely $W_m$; moreover, for every $i\le i_s$, we have $W_i=0$. Similarly, the injectivity of $\phi_0^{\ast}\cdots\phi_{i_{s+1}-1}^{\ast}$ implies that, for every $i\ge i_{s+1}$, we have $W_i=0$. It follows that 
$$
\sum_{i\in Q_0\setminus\{0\}}\theta_i\dim W_i=0,
$$
contradicting the $\bs\theta$-stableness of $(\mathbf{V},\bs\phi)$.
\end{proof}

\begin{Lem}\label{Lem-bar{V}bar{phi}}
Let $(\mathbf{V},\bs\phi)$ be a $\bs\theta$-stable representation representing a point of 
$$
\Rep^{\circ}(\bar{Q},\mathbf{d})^{\bs\theta-s}\cap\mu_{>0}^{-1}(\bs\zeta). 
$$
For any $0\le s\le \mu$, and any $i_s\le i<i_{s+1}$, define the subspace $W_i:=\phi^{\ast}_i\cdots\phi^{\ast}_{i_{s+1}-1}(V_{i_{s+1}})$. Then, the following assertions hold:
\begin{itemize}
\item[(i)] The maps $(\phi_i,\phi_i^{\ast})_{i_s\le i<i_{s+1}}$ induce maps $\bar{\bs\phi}:=(\bar{\phi}_i,\bar{\phi}_i^{\ast})_{i_s\le i<i_{s+1}}$ between the quotient spaces $\bar{\mathbf{V}}:=(\bar{V}_i:=V_i/W_i)_{i_s\le i<i_{s+1}}$, where $\bar{\phi}_{i_{s+1}-1}$ and $\bar{\phi}_{i_{s+1}-1}^{\ast}$ are the zero maps.
\item[(ii)] The pair $(\bar{\mathbf{V}},\bar{\bs\phi})$ forms a simple representation of the type A quiver with vertices $\{i_s,i_s+1,\ldots,i_{s+1}-1\}$, where $i_s$ is considered a framing vertex.
\item[(iii)] The representation $(\bar{\mathbf{V}},\bar{\bs\phi})$ satisfies
$$
(1+\bar{\phi}_{i-1}\bar{\phi}_{i-1}^{\ast})(1+\bar{\phi}_i^{\ast}\bar{\phi}_i)=\zeta_i
$$
for every $i_s<i<i_{s+1}$.
\end{itemize}
\end{Lem}
\begin{proof}
By definition, the maps $(\phi^{\ast}_i)_{i_s\le i<i_{s+1}}$ preserve the subspaces $W_i$ for $i_s\le i\le i_{s+1}$. The proof of equation (\ref{eq-Lem-ss->inj}) shows that $(\phi_i)_{i_s\le i<i_{s+1}}$ also preserve these subspaces. We see that the maps $(\bar{\phi}_i,\bar{\phi}_i^{\ast})_{i_s\le i<i_{s+1}}$ are well-defined. Now, we prove (ii). If $\bar{\mathbf{V}}'\subset\bar{\mathbf{V}}$ is a proper subrepresentation with $\bar{V}'_{i_s}=\bar{V}_{i_s}$, then there is a proper subrepresentation $\mathbf{V}'\subset\mathbf{V}$ which coincides with $\mathbf{V}$ away from the vertices $i_s<i<i_{s+1}$ and whose quotients recovers $\bar{\mathbf{V}}'$ at the vertices $i_s\le i<i_{s+1}$. For this subrepresentation, the equality in (\ref{eq-semi-st-2}) is attained, contradicting the stableness of $(\mathbf{V},\bs\phi)$. We then observe that there is no nontrivial subrepresentation $\bar{\mathbf{V}}'\subset\bar{\mathbf{V}}$ with $\bar{V}'_{i_s}=0$; this is due to the injectivity proved in Lemma \ref{Lem-st->inj}. Finally, (iii) follows from the corresponding equalities for the maps $(\phi_i,\phi_i^{\ast})_{i_s\le i<i_{s+1}}$.
\end{proof}

\begin{Lem}\label{Lem-M_j-V(jm)}
For any $0\le s\le \mu$, the endomorphism of $F_s/F_{s+1}$ induced by $\zeta_0(1+\phi^{\ast}_{0}\phi_{0})$ lies in $\bar{C}^{(s)}$, where $F_{\bullet}$, $V_{\bullet}$ and $\bs\phi$ are determined by $P$ and $\bar{C}$ in the manner described before Lemma \ref{Lem-Ker-subrepres}.
\end{Lem}
\begin{proof}
According to Lemma \ref{Lem-ss->inj}, the map $\phi_{0}^{\ast}\cdots\phi_{i_s-1}^{\ast}$ induces an isomorphism
$$
\iota_{s}:V_{i_s}\lisom F_s.
$$
We first show that under the isomorphism $\iota_{s}$, the restriction of $1+\phi^{\ast}_{0}\phi_{0}$ to $F_s$ is equal to $(\prod_{i=1}^{i_s}\zeta_{i})(1+\phi^{\ast}_{i_s}\phi_{i_s})$. We compute
\begingroup
\allowdisplaybreaks
\begin{align*}
(1+\phi^{\ast}_{0}\phi_{0})\phi^{\ast}_{0}\cdots\phi^{\ast}_{i_s-1}&=\phi^{\ast}_{0}(1+\phi_{0}\phi^{\ast}_{0})\phi^{\ast}_{1}\cdots\phi^{\ast}_{i_s-1}\\
&=\zeta_{1}\phi^{\ast}_{0}(1+\phi^{\ast}_{1}\phi_{1})\phi^{\ast}_{1}\cdots\phi^{\ast}_{i_s-1}.
\end{align*}
\endgroup
By induction, we have
$$
(1+\phi^{\ast}_{0}\phi_{0})\phi^{\ast}_{0}\cdots\phi^{\ast}_{i_s-1}=(\prod_{i=1}^{i_s}\zeta_{i})\phi^{\ast}_{0}\cdots\phi^{\ast}_{i_s-1}(1+\phi^{\ast}_{i_s}\phi_{i_s}),
$$
as asserted.

For any $i_{s-1}\le i<i_s$, we define 
$$
W_{i}=\Ima\phi^{\ast}_{i}\cdots\phi^{\ast}_{i_s-1}\subset V_{i}.
$$ 
The same proof of (\ref{eq-Lem-ss->inj}) shows that 
$$
\phi_{i}(W_{i})\subset W_{i+1},
$$ 
noting that the first induction step $\phi_{i_s-1}(W_{i_s-1})\subset W_{i_s}$ is automatic. It follows that $\bar{V}_{i}:=V_{i}/W_{i}$ and $(\bar{\phi}_{i},\bar{\phi}^{\ast}_{i})$ for $i_{s-1}\le i<i_s$, where $\bar{\phi}_{i}$ and $\bar{\phi}^{\ast}_{i}$ are the induced maps between quotient spaces, define an invertible representation of the double quiver of type $A$ with vertices $\{i\}_{i_{s-1}\le i<i_s}$. Since 
$$
(1+\bar{\phi}_{i-1}\bar{\phi}_{i-1}^{\ast})(1+\bar{\phi}_i^{\ast}\bar{\phi}_i)=\zeta_i
$$
for any $i_{s-1}<i< i_s$, it follows from \cite[Theorem 2.1]{CB04} that $\xi_{i_{s-1}}(1+\bar{\phi}^{\ast}_{i_{s-1}}\bar{\phi}_{i_{s-1}})$ lies in $\bar{C}^{(s-1)}$. It remains to note that $\xi_{i_{s-1}}=\zeta_0(\prod_{i=1}^{i_{s-1}}\zeta_{i})$.
\end{proof}

The above lemmas show that the morphism 
\begingroup
\allowdisplaybreaks
\begin{align}\label{eq-Lem-ss->inj-0}
\Rep^{\circ}(\bar{Q},\mathbf{d})^{\bs\theta-ss}\cap\mu_{>0}^{-1}(\bs\zeta)&\longrightarrow G_0\times G_0/P_0\\\nonumber
\bs\phi&\longmapsto (X=\zeta_0(\Id+\phi^{\ast}_{0}\phi_{0}),F_{\bullet}=(F_s)_{1\le s\le \mu}).
\end{align}
\endgroup
is well-defined and factors through $\mathbf{F}(\bar{C})$.

\begin{Prop}\label{Prop-F(C)-unframed-quiver}
The morphism (\ref{eq-Lem-ss->inj-0}) induces a $G_0$-equivariant isomorphism of algebraic varieties
$$
f:(\Rep^{\circ}(\bar{Q},\mathbf{d})^{\bs\theta-ss}\cap\mu_{>0}^{-1}(\bs\zeta))\ds G_{>0}\lisom\mathbf{F}(\bar{C}),
$$
where $G_{>0}:=\prod_{i\in Q_0\setminus\{0\}}\GL_{d_i}$.
\end{Prop}
\begin{proof}
The proposition will be proved by applying the Zariski main theorem, and this will be achieved in the following steps: (1) $\mathbf{F}(\bar{C})$ is a normal variety; (2) $(\Rep^{\circ}(\bar{Q},\mathbf{d})^{\bs\theta-ss}\cap\mu_{>0}^{-1}(\bs\zeta))\ds G_{>0}$ is connected; (3) $f$ is bijective.

We first show that $\mathbf{F}(\bar{C})$ is normal. The variety
$$
\{(X,gP_0)\in G_0\times G_0/P_0\mid X\in gP_0g^{-1}\}
$$
is a fibre bundle over $G_0/P_0$, with the fibre over $gP_0$ being the parabolic subgroup $gP_0g^{-1}$. Let $W\subset G_0/P_0$ be an open subset where this bundle trivialises; we may assume that the coset $\Id_V P_0\in W$. The trivial bundle $W\times P_0$ projects to $W\times L_0$ with fibre being $U_0$. It is easy to see that $\mathbf{F}(\bar{C})\cap W\times P_0$ is simply the inverse image of $\bar{C}\subset L_0$ under this projection. It follows from the normality of $\bar{C}$ that $\mathbf{F}(\bar{C})$ is normal. Besides, we see that $\mathbf{F}(\bar{C})$ is connected and the subset 
$$
\mathbf{F}(C):=\{(X,gP_0)\in G_0\times G_0/P_0\mid g^{-1}Xg\in CU_0\}
$$
is open dense in $\mathbf{F}(\bar{C})$. 

The connectedness of $(\Rep^{\circ}(\bar{Q},\mathbf{d})^{\bs\theta-ss}\cap\mu_{>0}^{-1}(\bs\zeta))\ds G_{>0}$ results from the connectedness of additive quiver varieties (see \cite[Proposition 3.21]{BS0}), as we explain in more detail below. Note that the multiplicative moment map $\mu_{>0}$ for a type $A$ quiver can be equivalently written as an additive moment map $\tilde{\mu}_{>0}$, as long as we modify the components of $\Rep^{\circ}(\bar{Q},\mathbf{d})$ by suitable scalars; i.e.,
$$
1+\phi_{i-1}\phi^{\ast}_{i-1}=\zeta_i(1+\phi^{\ast}_i\phi_i)
$$
is equivalent to
$$
\xi_{i-1}\phi_{i-1}\phi_{i-1}^{\ast}-\xi_i\phi^{\ast}_i\phi_i=\xi_i-\xi_{i-1}.
$$
This allows us to identify $\Rep^{\circ}(\bar{Q},\mathbf{d})^{\bs\theta-ss}\cap\mu_{>0}^{-1}(\bs\zeta)$ with an open subset of $\Rep(\bar{Q},\mathbf{d})^{\bs\theta-ss}\cap\tilde{\mu}_{>0}^{-1}(\bs\lambda)$, where $\bs\lambda$ is an additive deformation parameter determined by $\bs\xi$. Moreover, this open subset is defined as the nonvanishing locus of a function on $\Rep(\bar{Q},\mathbf{d})$ (i.e., the product of the determinants of all maps of the form $(1+\phi^{\ast}_i\phi_i)$). The level sets of a function are closed, and thus this function is constant on all orbit closures. It follows that $$\Rep^{\circ}(\bar{Q},\mathbf{d})^{\bs\theta-ss}\cap\mu_{>0}^{-1}(\bs\zeta)$$ is a saturated open subset, and thus we have the open embedding
$$
\Rep^{\circ}(\bar{Q},\mathbf{d})^{\bs\theta-ss}\cap\mu_{>0}^{-1}(\bs\zeta)\ds G_{>0}\hooklongrightarrow \Rep(\bar{Q},\mathbf{d})^{\bs\theta-ss}\cap\tilde{\mu}_{>0}^{-1}(\bs\lambda)\ds G_{>0}.
$$  
There is an isomorphism between the representation variety $\Rep(\bar{Q},\mathbf{d})$ and the representation variety for an unframed quiver with vertices $Q_0$, and the moment map of the latter is compatible with $\tilde{\mu}_{>0}$ under this isomorphism. This is known as Crawley-Boevey's trick \cite{CB01}. For the unframed quiver, define the deformation parameter at $0$ to be $-\sum_{i\in Q_0\setminus\{0\}}d_i\lambda_i$ and define the stability parameter at $0$ to be $-\sum_{i\in Q_0\setminus\{0\}}d_i\theta_i$. Then, we have an isomorphism of moment map fibres and compatible stability conditions; see \cite[\S 3.2 and Remark 5.1.4]{Ginz}. It follows that $\Rep(\bar{Q},\mathbf{d})^{\bs\theta-ss}\cap\tilde{\mu}_{>0}^{-1}(\bs\lambda)\ds G_{>0}$ is isomorphic to an additive quiver variety and is connected according to Bellamy-Schedler.

To prove the surjectivity of $f$, we show that every $(X,gP_0)\in\mathbf{F}(\bar{C})$ defines a quiver representation that is mapped to $(X,gP_0)$ by $f$; this is an analogue of \cite[Theorem 2.1]{CB04} for nontrivial stability conditions. The parabolic subgroup $gP_0g^{-1}$ determines a flag $F_{\bullet}$ as in (\ref{eq-subsec-Ch-var-multi-Q-1}). The choice of $(X,gP_0)$ with $X^{(s)}$ lying in $C^{(s)}$ for every $0\le s\le \mu$ defines a flag $V_{\bullet}$ refining $F_{\bullet}$ as in (\ref{eq-defn-V}), which amounts to a parabolic subgroup $\tilde{P}_0$ contained in $P_0$. Consider the proper morphism
$$
\pi:G_0\times G_0/\tilde{P}_0\longrightarrow G_0\times G_0/P_0
$$
which is the identity on the first component and the natural projection on the second component. Let $\bs\xi$ be the tuple of complex numbers determined by $\bar{C}$ as in the begining of this subsection, and let $\mathbf{F}(\bs\xi)$ be the closed subvariety of $G_0\times G_0/\tilde{P}_0$ consisting of $(X,g\tilde{P}_0)$ such that if $V_{\bullet}$ is the flag corresponding to $g\tilde{P}_0$, then $X$ induces the scalar endomorphism $\xi_i$ on $V_i/V_{i+1}$. It follows from \cite[Theorem 2.1]{CB04} that the endomorphism $X^{(s)}$ that $X$ induces on $F_s/F_{s+1}$ lies in $\bar{C}^{(s)}$, which shows that $\pi$ maps $\mathbf{F}(\bs\xi)$ to $\mathbf{F}(\bar{C})$. We have seen that if $(X,gP_0)$ lies in $\mathbf{F}(C)$, then it lifts to an element of $\mathbf{F}(\bs\xi)$. By the properness of $\pi$, we deduce that $\pi(\mathbf{F}(\bs\xi))=\mathbf{F}(\bar{C})$; that is, every element of $\mathbf{F}(\bar{C})$ lifts to some $(X,g\tilde{P}_0)$. We then interpret $(X,g\tilde{P}_0)$ as a quiver representation. Regard $g\tilde{P}_0$ as a flag $V'_{\bullet}$ on $V'_0=\mathbb{C}^n$, then $\dim V'_i=d_i$. For $0\le i\le \nu-1$, define $\phi_i^{\ast}:V'_{i+1}\rightarrow V'_i$ to be the inclusion and define $\phi_i:V'_i\rightarrow V'_{i+1}$ as $(\xi_i^{-1}X-\Id)$. It is easy to see that $(V'_{\bullet},\bs\phi=(\phi_i,\phi_i^{\ast})_{0\le i\le \nu-1})$ is a well-defined quiver representation and lies in $\mu_{>0}^{-1}(\bs\zeta)$. It is also easy to verify that $(V'_{\bullet},\bs\phi)$ is $\bs\theta$-semi-stable. Indeed, there is no nonzero subrepresentation $W$ with $W_0=0$, since all $\phi^{\ast}_i$ are injective, and for every subrepresentation $W$ with $W_0=V_0$, the inequality
$$
\sum_{i\in Q_0\setminus\{0\}}\theta_i\dim W_i\le\sum_{i\in Q_0\setminus\{0\}}d_i\theta_i
$$
is automatically true since all $\theta_i$ are nonnegative. We have thus proved that $f$ is surjective.

As our last step, we prove that $f$ is injective. Recall that a $\bs\theta$-polystable representation $(\mathbf{V},\bs\phi)$ admits a decomposition of the form
$$
(\mathbf{V},\bs\phi)=(\mathbf{V}^{(1)},\bs\phi^{(1)})\oplus(\mathbf{V}^{(2)},\bs\phi^{(2)}),
$$
where $(\mathbf{V}^{(1)},\bs\phi^{(1)})$ is $\bs\theta$-stable with $V^{(1)}_0=V_0$ and $(\mathbf{V}^{(2)},\bs\phi^{(2)})$ is $\bs\theta$-polystable with $V^{(2)}_0=0$. Let us begin by showing that if $(X,gP_0)=f((\mathbf{V},\bs\phi))$, then $\mathbf{d}^{(1)}$, the dimension vector of $(\mathbf{V}^{(1)},\bs\phi^{(1)})$, is determined by $(X,gP_0)$.

We regard $gP_0$ as a flag $(F_s)_{0\le s\le\mu}$ on $V_0$. By the definition of $f$, we have $\dim F_s=\dim V_{i_s}$. The $\bs\theta$-semistableness of $(\mathbf{V},\bs\phi)$ implies that the components of $\mathbf{V}^{(2)}$ at any $i_s$ is zero. We deduce that $\dim V^{(1)}_{i_s}=\dim F_s$ for any $0<s\le \mu$. It remains to show that for any $0\le s\le \mu$ and any $i_s<i<i_{s+1}$, the dimension of $V^{(1)}_i$ is determined by $(X,gP_0)$. For any $s$, we apply Lemma \ref{Lem-bar{V}bar{phi}} to the $\bs\theta$-stable $(\mathbf{V}^{(1)},\bs\phi^{(1)})$ and obtain a simple representation of the type A quiver with vertices $\{i_s,i_s+1,\ldots,i_{s+1}-1\}$. According to the proof of Lemma \ref{Lem-M_j-V(jm)}, the map $1+\bar{\phi}^{\ast}_{i_s}\bar{\phi}_{i_s}$ is equal to the map induced by $1+\phi_0^{\ast}\phi_0$ on $F_s/F_{s+1}$ up to a scalar, but the latter map is equal to $X$ up to a scalar; therefore, $1+\bar{\phi}^{\ast}_{i_s}\bar{\phi}_{i_s}$ is determined by $(X,gP_0)$. It is known that the dimension vector of simple representation of a type A quiver is determined by $\bs\zeta$ and $1+\bar{\phi}^{\ast}_{i_s}\bar{\phi}_{i_s}$. This proves our claim that $\mathbf{d}^{(1)}$ is determined by $(X,gP_0)$.

We now complete the proof of the injectivity of $f$. The dimension vector of $\mathbf{V}^{(2)}$ is determined by $\mathbf{d}$ and $\mathbf{d}^{(1)}$. The isomorphism class of $\mathbf{V}^{(2)}$ is then determined by its dimension vector, since it is supported on a union of type A unframed quivers. The isomorphism class of $\mathbf{V}^{(1)}$ can be determined as follows. By Lemma \ref{Lem-ss->inj} and Lemma \ref{Lem-st->inj}, the map $\phi_{0}^{\ast}\cdots\phi_{i-1}^{\ast}$ is injective for every $0<i\le\nu$. Up to the conjugation by a suitable element of $\prod_{0<i\le\nu}\GL_{d_i}$, we map assume that $\phi^{\ast}_i$ is the inclusion map for every $i$. With the prescribed $X$, there is a unique map $\phi_0$ satisfying $X=\xi_0(1+\phi^{\ast}_0\phi_0)$. Inductively, this determines all $\phi_i$ and thus the isomorphism class of $\mathbf{V}^{(1)}$.
\end{proof}

\subsection{Multiplicative quiver varieties for comet-shaped quivers}\label{subsec-multi-comet}\hfill

In the rest of this article, we will only consider multiplicative quiver varieties for the comet-shaped quiver below.
\[\begin{tikzcd}
	& {[1,1]} & {[1,2]} & \cdots & {[1,\nu_1]} \\
	& {[2,1]} & {[2,2]} & \cdots & {[2,\nu_2]} \\
	\star & \vdots & \vdots && \vdots \\
	& {[k,1]} & {[k,2]} & \cdots & {[k,\nu_k]}
	\arrow[from=1-2, to=1-3]
	\arrow[from=1-3, to=1-4]
	\arrow[from=1-4, to=1-5]
	\arrow[from=2-2, to=2-3]
	\arrow[from=2-3, to=2-4]
	\arrow[from=2-4, to=2-5]
	\arrow[from=3-1, to=1-2]
	\arrow[from=3-1, to=2-2]
	\arrow[shift left, from=3-1, to=3-1, loop, in=150, out=210, distance=10mm]
	\arrow["\cdots"', shift right, from=3-1, to=3-1, loop, in=145, out=215, distance=20mm]
	\arrow[from=3-1, to=4-2]
	\arrow[from=4-2, to=4-3]
	\arrow[from=4-3, to=4-4]
	\arrow[from=4-4, to=4-5]
\end{tikzcd}\]
The vertex set $Q_0$ is the union of a central vertex $\star$ and $\{[j,i]\}_{1\le j\le k,~1\le i\le \nu_j}$, and we will sometimes represent $\star$ by $[j,0]$ for some $1\le j\le k$; the set of arrows $Q_1$ consists of an arrow $[j,i]\rightarrow[j,i+1]$ for every $1\le j\le k$ and every $0\le i\le \nu_j-1$. The loop-free subquiver consisting of the vertices $\{[j,i]\}_{0\le i\le \nu_j}$ and the arrows between them is called a leg of the quiver. The components of a representation $(\mathbf{V},\bs\phi)$ of the double $\bar{Q}$ will be written as $\phi_{[j,i]}:V_{[j,i]}\rightarrow V_{[j,i+1]}$ and $\phi_{[j,i]}^{\ast}:V_{[j,i+1]}\rightarrow V_{[j,i]}$. 

Now, we introduce the moduli space $\mathcal{M}_B(\mathbf{P},\bar{\mathcal{C}},\bs\beta)$ as in the statement of our main theorem and explain how it is identified with an open subset of a multiplicative quiver variety. The moduli space in question parametrises filtered local systems, for which we introduce some notation. As in the introduction, we denote by $X=\bar{X}\setminus\{x_1,\ldots,x_k\}$ a Riemann surface of genus $g$ with $k$ punctures. Denote by $\Pi=\pi_1(X,x_0)$ its fundamental group with a base point $x_0$. For each $1\le j\le k$, let $U_j\subset\bar{X}$ be a simply connected open neighbourhood of $x_j$ and write $U_j^{\ast}=U_j\setminus\{x_j\}$. Let $y_j\in U_j^{\ast}$ and $\gamma'_j$ a loop in $U_j^{\ast}$ based at $y_j$ which generates $\pi_1(U_j^{\ast})$. Then, we may choose $\{\alpha_i,\beta_i\}_{1\le i\le g}\subset\Pi$ and paths $l_j$ connecting $x_0$ and $y_j$ so that $\Pi$ is generated by $\{\alpha_i,\beta_i,\gamma_j:=l_j^{-1}\gamma'_jl_j\mid 1\le i\le g, 1\le j\le k\}$ subject to the unique relation $\prod_i[\alpha_i,\beta_i]\prod_j\gamma_j=1$.

A filtered local system, as was recalled in \S \ref{subsec-main-thm}, is equivalent to a filtered representation of $\Pi$; that is, a homomorphism $\rho:\Pi\rightarrow\GL(V)$ equipped with $k$ flags on $V$, each denoted by
\begin{equation}\label{eq-flag-comet-quiver}
F_{[j,\bullet]}=(V=F_{[j,0]}\supset F_{[j,1]}\supset\cdots\supset F_{[j,\mu_j]}\supset F_{[j,\mu_j+1]}=0),
\end{equation}
such that each $\rho(\gamma_j)$ preserves $F_{[j,\bullet]}$. We will denote by $(\rho,F_{[\bullet,\bullet]})$ a filtered representation of $\Pi$.  Since $\rho(\gamma_j)$ preserves $F_{[\bullet,\bullet]}$, it determines an element $\bar{M}_j\in\GL(\Gr F_{[j,\bullet]})$, called the graded monodromy at $x_j$.

A nice moduli space is obtained if we fix the types of the flags and impose suitable constraints on the monodromy. The subgroup of $\GL(V)$ preserving a flag $F_{[j,\bullet]}$ is a parabolic subgroup $P_j\subset\GL(V)$, and fixing the type of this flag means fixing the conjugacy class of this parabolic subgroup. The parabolic subgroups conjugate to $P_j$ are parametrised by $\GL(V)/P_j$. Let $L_j$ be the Levi quotient of $P_j$ and $U_j$ the unipotent radical of $P_j$. We will write $\mathbf{P}=(P_j)_{1\le j\le k}$. For each $j$, let $\bar{C}_j$ be the closure of an $L_j$-conjugacy class $C_j\subset L_j$, and denote by $\bar{\mathcal{C}}=(\bar{C}_j)_{1\le j\le k}$ the tuple of these closures of conjugacy classes. 

Every point of the following algebraic variety represents a filtered representation
\begin{equation}\label{eq-rep-var-filter}
\Rep_B(g,\mathbf{P},\bar{\mathcal{C}}):=\{((A_i,B_i)_i,(X_j,g_jP_j)_j)\in \GL_n^{2g}\times\prod_{j=1}^k\mathbf{F}(\bar{C}_j)\mid\prod_{i=1}^g[A_i,B_i]\prod_{j=1}^kX_j=\Id\},
\end{equation}
where $\mathbf{F}(\bar{C}_j)$ is defined by $P_j$ and $\bar{C}_j$ as in (\ref{eq-defn-F(C)}). Two representations are isomorphic if and only if the corresponding points are $\GL(V)$-conjugate. We introduce a polarisation as follows. For any $1\le j\le k$ and $0\le s\le \mu_j$, we write $f_{[j,s]}=\dim F_{[j,s]}$. Let $\bs\beta=(\beta_{[j,s]})_{1\le j\le k,0\le s\le \mu_j}$ be a tuple of rational numbers satisfying $\beta_{[j,s]}<\beta_{[j,s']}$ whenever $s<s'$ and 
\begin{equation}\label{eq-beta-deg=0}
\sum_{j=1}^k\sum_{s=0}^{\mu_j}\beta_{[j,s]}(f_{[j,s]}-f_{[j,s+1]})=0.
\end{equation}
Let $K$ be the smallest positive integer such that $b_{[j,s]}:=\beta_{[j,s]}K$ is an integer for all $[j,s]$. We then define $c_{[j,s]}:=b_{[j,s]}-b_{[j,s-1]}$ for $s>0$. For every $[j,s]$ with $s>0$, there is a projection from $\Rep_B(g,\mathbf{P},\bar{\mathcal{C}})$ to the Grassmannian $\Gr(V,f_{[j,s]})$; the very ample line bundle $\mathcal{O}(c_{[j,s]})$ on $\Gr(V,f_{[j,s]})$ pullbacks to a line bundle on $\Rep_B(g,\mathbf{P},\bar{\mathcal{C}})$, and the tensor product of these line bundles for all $[j,s]$ with $s>0$ is a $\GL(V)$-equivariant line bundle $\mathcal{L}$, which is the desired polarisation. The Hilbert-Mumford criterion then translates the $\mathcal{L}$-stability into a numerical one: the computation needed only involves the flag variety components of $\Rep_B(g,\mathbf{P},\bar{\mathcal{C}})$, which turns out to be identical to what one finds in similar contexts like parabolic vector bundles, and the result is precisely the $\bs\beta$-stability condition (\ref{eq-filt-locsys-stability}) defined as in the introduction. The moduli space of filtered local system is defined to be the GIT-quotient of $\Rep_B(g,\mathbf{P},\bar{\mathcal{C}})$ with the polarisation $\mathcal{L}$:
\begin{equation}\label{eq-char-var-filter}
\mathcal{M}_B(g,\mathbf{P},\bar{\mathcal{C}},\bs\beta):=\Rep_B(g,\mathbf{P},\bar{\mathcal{C}})~\ds_{\mathcal{L}}\GL(V).
\end{equation}
This quotient parametrises the closed orbits in the $\bs\beta$-semi-stable locus $\Rep_B(g,\mathbf{P},\bar{\mathcal{C}})^{\bs\beta-ss}$. We will always work with a fixed Riemann surface and consider the moduli spaces associated to varying invariants $\mathbf{P}$, $\bar{\mathcal{C}}$, $\bs\beta$, and thus $g$ will usually be omitted from the notation.

We will show that $\mathcal{M}_B(\mathbf{P},\bar{\mathcal{C}},\bs\beta)$ is an open subvariety of a multiplicative quiver variety $\mathcal{M}(\bs\zeta,\mathbf{d},\bs\theta)$. The quiver data $\bs\zeta$, $\mathbf{d}$ and $\bs\theta$ are determined by $\mathbf{P}$, $\bar{\mathcal{C}}$, $\bs\beta$ as follows. According to \S \ref{subsec-multi-Q-typeA}, for each $1\le j\le k$, there exists a type $A$ quiver with vertices $\{[j,i]\}_{0\le i\le \nu_j}$, a tuple of complex numbers $\bs\zeta_j=(\zeta_{[j,i]})_{0\le i\le \nu_j}$ and a dimension vector $\mathbf{d}_j=(d_{[j,i]})_{0\le i\le \nu_j}$, so that the variety $\mathbf{F}(\bar{C}_j)$ is isomorphic to a multiplicative quiver variety; see Proposition \ref{Prop-F(C)-unframed-quiver}. Identifying the vertices $\{[j,0]\}_{1\le j\le k}$ of the resulting type $A$ quivers, we obtain a star-shaped quiver with $k$ legs; we then define a comet-shaped quiver $Q$ by joining $g$ loops to the vertex $\star$. Since $d_{[j,0]}=n$ for all $j$, the dimension vectors $\mathbf{d}_j$ glue and forms a dimension vector $\mathbf{d}$ supported on $Q$. We then define $\bs\zeta$ in such a way that it agrees with $\bs\zeta_j$ on the $j$-th leg and $\zeta_{\star}=\prod_{j=1}^k\zeta_{[j,0]}$. Note that Proposition \ref{Prop-F(C)-unframed-quiver} only requires $\bs\theta=(\theta_i)_{i\in Q_0}$ to satisfy $\theta_i\in\mathbb{Q}_{>0}$ whenever $i$ lies in a particular subset $\{i_s\}_{1\le s\le \mu_j}\subset\{1,2,\ldots,\nu_j\}$ and $\theta_i=0$ otherwise. We now further require that $\bs\theta$ is determined by $\bs\beta$ in the following manner. Recall that the conjugacy class of $P_j$ determines the length $\mu_j$ of the flag (\ref{eq-flag-comet-quiver}), and $\bs\beta$ assigns to $F_{[j,s]}$ the number $\beta_{[j,s]}$. We define $\bs\theta$ so that 
\begin{itemize}
\item $\theta_{[j,i_s]}=\beta_{[j,s]}-\beta_{[j,s-1]}$ for every $j$ and every $1\le s\le \mu_j$,
\item $\theta_{[j,i]}=0$ if $i\ne\star$ and $i\ne i_s$ for any $1\le s\le \mu_j$, and
\item $\theta_{\star}d_{\star}+\sum_{j=1}^k\sum_{1\le i\le \nu_j}\theta_{[j,i]}d_{[j,i]}=0$.
\end{itemize}
In view of (\ref{eq-beta-deg=0}), the last equality is equivalent to $\theta_{\star}=\sum_{1\le j\le k}\beta_{[j,0]}$.

\begin{Prop}\label{Prop-B-open-multi}
Let $Q$, $\bs\zeta$, $\mathbf{d}$ and $\bs\theta$ be determined by $g$, $\mathbf{P}$, $\bar{\mathcal{C}}$ and $\bs\beta$ as above. Then, there is a natural open embedding
$$
\mathcal{M}_B(\mathbf{P},\bar{\mathcal{C}},\bs\beta)\hooklongrightarrow\mathcal{M}(\bs\zeta,\mathbf{d},\bs\theta),
$$
which is an isomorphism if $g=0$.
\end{Prop}
\begin{proof}
We first introduce the following locally closed subvariety
$$
\Rep^{\circ}(g,\mathbf{P},\bar{\mathcal{C}})\subset \gl_n^{2g}\times\prod_{j=1}^k\mathbf{F}(\bar{C}_j)
$$ 
consisting of the tuples $((\psi_i,\psi_i^{\ast})_{1\le i\le g},(X_j,gP_j)_{1\le j\le k})$ satisfying the following conditions:
\begin{equation}\label{eq-Prop-B-open-multi-01}
(1+\psi_i\psi_i^{\ast})\text{ and }(1+\psi_i^{\ast}\psi_i)\text{ are invertible for any }i;
\end{equation}
and
\begin{equation}\label{eq-Prop-B-open-multi-02}
\prod_{i=1}^g(1+\psi_i\psi_i^{\ast})(1+\psi_i^{\ast}\psi_i)^{-1}\prod_{j=1}^kX_j=1.
\end{equation}
Define an open embedding
\begin{equation}
\iota:\Rep_B(g,\mathbf{P},\bar{\mathcal{C}})\longrightarrow \Rep^{\circ}(g,\mathbf{P},\bar{\mathcal{C}}).
\end{equation}
as follows. Let $((A_i,B_i)_i,(X_j,g_jP_j)_j)\in\Rep_B(g,\mathbf{P},\bar{\mathcal{C}})$. Its image under $\iota$ has the same components $X_j$ and $g_jP_j$ for every $j$. For any $1\le i\le g$, define 
\begin{equation}\label{eq-Prop-B-open-multi-2}
\psi_i=B_i^{-1},\text{ and }\psi_i^{\ast}=B_iA_i-B_i,
\end{equation}
so that
\begin{equation}\label{eq-Prop-B-open-multi-3}
A_i=1+\psi_i\psi_i^{\ast},\text{ and }B_i=\psi_i^{-1}.
\end{equation}
We have
$$
(1+\psi_i\psi_i^{\ast})(1+\psi_i^{\ast}\psi_i)^{-1}=A_iB_iA_i^{-1}B_i^{-1},
$$
and thus $\iota$ is well-defined. We claim that $\iota$ induces an isomorphism onto an open subset $U\subset\Rep^{\circ}(g,\mathbf{P},\bar{\mathcal{C}})$, where
\begin{equation}\label{eq-Prop-B-open-multi-U}
U:=\{((\psi_i,\psi_i^{\ast})_{1\le i\le g},(X_j,gP_j)_{1\le j\le k})\in\Rep^{\circ}(g,\mathbf{P},\bar{\mathcal{C}})\mid \psi_i\text{ is invertible for any }i\}.
\end{equation}
Indeed, we may define a morphism $U\rightarrow\Rep_B(g,\mathbf{P},\bar{\mathcal{C}})$ by the formula (\ref{eq-Prop-B-open-multi-3}), which is obviously inverse to $\iota$. 

The open embedding $\iota$ induces an isomorphism 
$$
\Rep_B(g,\mathbf{P},\bar{\mathcal{C}})^{\bs\beta-ss}\lisom U\cap\Rep^{\circ}(g,\mathbf{P},\bar{\mathcal{C}})^{\bs\beta-ss};
$$
note that $(\psi_i,\psi_i^{\ast})_{1\le i\le g}$ preserves a subspace $W\subset\mathbb{C}^{n}$ if and only if $(A_i,B_i)_{1\le i\le g}$ does under the correspondences (\ref{eq-Prop-B-open-multi-2}) and (\ref{eq-Prop-B-open-multi-3}). In order to show that this isomorphism induces an open embedding
$$
\mathcal{M}_B(\mathbf{P},\bar{\mathcal{C}},\bs\beta)\hookrightarrow \Rep^{\circ}(g,\mathbf{P},\bar{\mathcal{C}})^{\bs\beta-ss}\ds \GL_n,
$$
we need to show that a closed orbit in $U$ is also closed in $\Rep^{\circ}(g,\mathbf{P},\bar{\mathcal{C}})^{\bs\beta-ss}$. Our arguments are analogous to the proof of \cite[Lemma 5.8]{Boa15}. The open subset $U$ is defined by the nonvanishing of the $\GL_n$-invariant function $f:=\prod_{i=1}^g\det(\phi_i)$ on $\Rep^{\circ}(g,\mathbf{P},\bar{\mathcal{C}})$. If $f$ is nonzero (necessarily of constant value) on an orbit $\mathcal{O}\subset U$, then it takes the same value on its closure $\overline{\mathcal{O}}$, and thus $\overline{\mathcal{O}}$ is contained in $U$; that is, if $\mathcal{O}$ is not closed in $\Rep^{\circ}(g,\mathbf{P},\bar{\mathcal{C}})^{\bs\beta-ss}$, then it is not closed in $U$.

The last step of the proof is to show that there is a $\GL_n$-equivariant isomorphism
\begin{equation}\label{eq-Prop-B-open-multi-4}
\Rep^{\circ}(g,\mathbf{P},\bar{\mathcal{C}})^{\bs\beta-ss}\lisom\Rep^{\circ}(\bar{Q},\bs\zeta,\mathbf{d})^{\bs\theta-ss}\ds G_{>0},
\end{equation}
where the $\bs\beta$-semi-stable locus of $\Rep^{\circ}(g,\mathbf{P},\bar{\mathcal{C}})$ is defined in the same way as filtered local systems, and $G_{>0}:=\prod_{1\le j\le k,~i>0}\GL_{d_{[j,i]}}$. Since the quotient of $\Rep^{\circ}(\bar{Q},\bs\zeta,\mathbf{d})^{\bs\theta-ss}\ds G_{>0}$ by $\GL_n$ is isomorphic to $\mathcal{M}(\bs\zeta,\mathbf{d},\bs\theta)$, we obtain the desired open embedding.

To prove the isomorphism (\ref{eq-Prop-B-open-multi-4}), we proceed as follows. For $1\le j\le k$, let $Q^{(j)}$ be a leg of $Q$, which is a loop-free quiver with the set of vertices $Q^{(j)}_0=\{[j,i]\mid 0\le i\le \nu_j\}$ (recall that $[j,0]$ is simply the vertex $\star$). Denote by $\mathbf{d}^{(j)}$ the restriction of $\mathbf{d}$ to $Q^{(j)}$ and by $\bs\zeta^{(j)}$ and $\bs\theta^{(j)}$ the respective restrictions of $\bs\zeta$ and $\mathbf{d}$ to $Q^{(j)}_0\setminus\{\star\}$. Applying Proposition \ref{Prop-F(C)-unframed-quiver} to every leg of $Q$, we obtain an isomorphism
\begin{equation}\label{eq-Prop-B-open-multi-5}
\gl_n^{2g}\times\prod_{j=1}^k(\Rep^{\circ}(\bar{Q}^{(j)},\mathbf{d}^{(j)})^{\bs\theta^{(j)}-ss}\cap\mu_{>0}^{-1}(\bs\zeta^{(j)}))\ds G^{(j)}_{>0}\lisom\gl_n^{2g}\times\prod_{j=1}^n\mathbf{F}(\bar{C}_j).
\end{equation}
Denote by 
$$
\Rep^{\circ}(\bar{Q},\bs\zeta,\mathbf{d})^{(\bs\theta^{(1)},\ldots,\bs\theta^{(k)})-ss}\subset \gl_n^{2g}\times\prod_{j=1}^k\Rep^{\circ}(\bar{Q}^{(j)},\mathbf{d}^{(j)})^{\bs\theta^{(j)}-ss}\cap\mu_{>0}^{-1}(\bs\zeta^{(j)})
$$  
the locally closed subset defined by the invertibility of the $(\psi_i,\psi_i^{\ast})$-components and the closed condition $\mu_{\star}((\mathbf{V},\bs\phi))=\zeta_{\star}$. The isomorphism (\ref{eq-Prop-B-open-multi-5}) then restricts to
$$
\Rep^{\circ}(\bar{Q},\bs\zeta,\mathbf{d})^{(\bs\theta^{(1)},\ldots,\bs\theta^{(k)})-ss}\ds G_{>0}\lisom\Rep^{\circ}(g,\mathbf{P},\bar{\mathcal{C}}).
$$
We will show that $\Rep^{\circ}(\bar{Q},\bs\zeta,\mathbf{d})^{\bs\theta-ss}$ is precisely the inverse image of $\Rep^{\circ}(g,\mathbf{P},\bar{\mathcal{C}})^{\bs\beta-ss}$ under the quotient by $G_{>0}$, from which the isomorphism (\ref{eq-Prop-B-open-multi-4}) follows.

It is easy to see that $\Rep^{\circ}(\bar{Q},\bs\zeta,\mathbf{d})^{\bs\theta-ss}$ is an open subset of $\Rep^{\circ}(\bar{Q},\bs\zeta,\mathbf{d})^{(\bs\theta^{(1)},\ldots,\bs\theta^{(k)})-ss}$; indeed, for $(\mathbf{V},\bs\phi)\in\Rep^{\circ}(\bar{Q},\bs\zeta,\mathbf{d})^{\bs\theta-ss}$, its restriction to $Q^{(j)}$, regarded as a quiver with framing vertex $\star$, is a $\bs\theta^{(j)}$-semi-stable representation. Let us show that $\Rep^{\circ}(\bar{Q},\bs\zeta,\mathbf{d})^{\bs\theta-ss}$ maps into $\Rep^{\circ}(g,\mathbf{P},\bar{\mathcal{C}})^{\bs\beta-ss}$. Let $(\mathbf{V},\bs\phi)$ be a $\bs\theta$-semi-stable representation as above. Let $W_{\star}\subset V_{\star}$ be a subspace that is preserved by $(\psi_i,\psi_i^{\ast})_{1\le i\le g}$ and $(X_j)_{1\le j\le k}$. We define a subrepresentation $\mathbf{W}\subset\mathbf{V}$ by 
\begin{equation}\label{eq-defn-Wstar-repres}
W_{[j,i]}=(\phi_{[j,0]}^{\ast}\phi_{[j,1]}^{\ast}\cdots\phi_{[j,i-1]}^{\ast})^{-1}(W_{\star}).
\end{equation}
For each $j$, let $\{i_s\}_{1\le s\le \mu_j}\subset\{1,2,\ldots,\nu_j\}$ be the subset such that $\{[j,i_s]\}_{1\le j\le k,~1\le s\le \mu_j}$ are precisely the vertices with $\theta_{[j,i]}>0$. Then, the $\bs\theta$-semi-stability gives
$$
\theta_{\star}\dim W_{\star}+\sum_{j=1}^k\sum_{s=1}^{\mu_j}\theta_{[j,i_s]}\dim W_{[j,i_s]}\le 0.
$$
For any $[j,i_s]$, we have $W_{[j,i_s]}\cong W_{\star}\cap F_{[j,s]}$ by Lemma \ref{Lem-ss->inj} (recall that $F_{[j,s]}$ is the image of $V_{[j,i_s]}$ under $\phi_{[j,0]}^{\ast}\phi_{[j,1]}^{\ast}\cdots\phi_{[j,i_s-1]}^{\ast}$). It follows from the definition of $\bs\theta$ that the above inequality is equivalent to the $\bs\beta$-semi-stability of the filtered representation with underlying vector space $V_{\star}$. We then verify that if a quiver representation $(\mathbf{V},\bs\phi)$ in $\Rep^{\circ}(\bar{Q},\bs\zeta,\mathbf{d})^{(\bs\theta^{(1)},\ldots,\bs\theta^{(k)})-ss}$ is sent to a $\bs\beta$-semi-stable filtered representation, then $(\mathbf{V},\bs\phi)$ is $\bs\theta$-semi-stable. Let $\mathbf{W}\subset\mathbf{V}$ be a subrepresentation. Let $\tilde{\mathbf{W}}$ be the largest subrepresentation with $\tilde{W}_{\star}=W_{\star}$, which is defined by (\ref{eq-defn-Wstar-repres}). We have $\mathbf{W}\subset \tilde{\mathbf{W}}$. The same reasoning as above shows that the $\bs\beta$-semi-stability of filtered representation gives
$$
\theta_{\star}\dim \tilde{W}{\star}+\sum_{j=1}^k\sum_{s=1}^{\mu_j}\theta_{[j,i_s]}\dim \tilde{W}_{[j,i_s]}\le 0.
$$
Since $\dim W_{[j,i_s]}\le\dim \tilde{W}_{[j,i_s]}$ for any $[j,i_s]$, we obtain the desired inequality for $\mathbf{W}$.
\end{proof}

\begin{Cor}
The character variety $\mathcal{M}_B(\mathbf{P},\bar{\mathcal{C}},\bs\beta)$ is normal.
\end{Cor}
\begin{proof}
This follows from \cite{CB03} and \cite[Theorem 5.4]{KS}.
\end{proof}
\begin{Cor}\label{Cor-dense-M_B^s}
The locus of $\bs\beta$-stable filtered local systems $\mathcal{M}_B(\mathbf{P},\bar{\mathcal{C}},\bs\beta)^s$ is dense in $\mathcal{M}_B(\mathbf{P},\bar{\mathcal{C}},\bs\beta)$ if $\mathbf{d}\in\Sigma_{\bs\zeta,\bs\theta}$, where $\bs\zeta$, $\mathbf{d}$ and $\bs\theta$ are as in Proposition \ref{Prop-B-open-multi}.
\end{Cor}
\begin{proof}
By \cite[Proposition 2.19 and Remark 2.20]{ST}, the stable locus $\mathcal{M}(\bs\zeta,\mathbf{d},\bs\theta)^s$ is dense in $\mathcal{M}(\bs\zeta,\mathbf{d},\bs\theta)$, from which the assertion follows.
\end{proof}

\section{Point-counting for almost generic character varieties}\label{sec-RGSSM}

In this section, we will prove the irreducibility of affine character varieties with almost generic semi-simple monodromies $\mathcal{C}=(C_j)_{1\le j\le k}$ by counting points over finite fields. The representation variety (\ref{eq-rep-var-filter}) becomes
\begin{equation}\label{eq-sec-RGSSM-1}
\Rep_B(\mathcal{C}):=\{((A_i,B_i)_i,(X_j)_j)\in \GL_n^{2g}\times\prod_{j=1}^kC_j\mid\prod_{i=1}^g[A_i,B_i]\prod_{j=1}^kX_j=\Id\},
\end{equation}
while the character variety (\ref{eq-char-var-filter}) is now the affine GIT quotient 
\begin{equation}\label{eq-sec-RGSSM-2}
\mathcal{M}_B(\mathcal{C}):=\Rep_B(\mathcal{C})\ds\GL_n.
\end{equation}
The quotient stack
\begin{equation}\label{eq-sec-RGSSM-3}
\mathfrak{M}_B(\mathcal{C}):=[\Rep_B(\mathcal{C})/\GL_n]
\end{equation}
will be called the character stack associated to $g$ and $\mathcal{C}$. The key ingredient for our point-counting problem is Scognamiglio's formula for non-generic character stacks, which we will review in \S \ref{subsec-Scog}. 

\subsection{Plethysm}\label{subsec-pleth}\hfill

The statement of Scognamiglio's formula involves the plethystic exponential, which we recall below. A lambda ring is a commutative ring $R$ over $\mathbb{Q}$ equipped with $\mathbb{Q}$-linear maps $p_d:R\rightarrow R$ for all $d\in\mathbb{Z}_{\ge1}$ that satisfy certain properties. Given a lambda ring $R$, consider $R[[T_i]]_{i\in I}$, the ring of formal power series in a finite set of independent variables indexed by a finite set $I$ with coefficients in $R$. For any $\mathbf{d}\in\mathbb{Z}_{\ge0}^I$, write $T^{\mathbf{d}}=\prod_{i\in I}T_i^{d_i}$. We extend the lambda ring structure on $R$ to $R[[T_i]]_{i\in I}$ by defining $p_d(rT^{\mathbf{d}})=p_d(r)T^{d\mathbf{d}}$ for any $r\in R$ and any $\mathbf{d}\in\mathbb{Z}_{\ge0}^I$. We will primarily consider $R=\mathbb{Q}(q)$, the field of rational functions in $q$, and equip it with a lambda ring structure by defining $p_d(f(q))=f(q^d)$. 

Denote by $R[[T_i]]^+_{i\in I}\subset R[[T_i]]_{i\in I}$ the ideal generated by $\{T_i\}_{i\in I}$. The plethystic exponential is the following map
\begingroup
\allowdisplaybreaks
\begin{align*}
\Exp:R[[T_i]]^+_{i\in I}&\longrightarrow 1+ R[[T_i]]^+_{i\in I}\\
F&\longmapsto \exp\big(\sum_{n\ge1}\frac{p_n(F)}{n}\big),
\end{align*}
\endgroup
where the usual exponential $\exp$ is
\begingroup
\allowdisplaybreaks
\begin{align*}
\exp:R[[T_i]]^+_{i\in I}&\longrightarrow 1+ R[[T_i]]^+_{i\in I}\\
f&\longmapsto 1+\sum_{n\ge1}\frac{f^n}{n!}.
\end{align*}
\endgroup
For any $n\in\mathbb{Z}_{\ge0}$, denote by $\mathcal{P}(n)$ the set of partitions of $n$ (so that $\mathcal{P}(0)$ consists of the empty partition). Define
$$
\sigma_n(F)=\sum_{\lambda\in\mathcal{P}(n)}\frac{p_{\lambda}(F)}{z_{\lambda}},
$$
for any $F\in R[[T_i]]_{i\in I}$, where
$$
p_{\lambda}(F):=p_{\lambda_1}(F)\cdots p_{\lambda_l}(F),
$$
for $\lambda=(\lambda_1,\ldots,\lambda_l)$, and
$$
z_{\lambda}:=\prod_jm_j!j^{m_j},
$$
if we write $\lambda=(1^{m_1},2^{m_2},\ldots)$.

\subsection{Scognamiglio's formula}\label{subsec-Scog}\hfill

The basis of the study of character varieties via point-counting is the fact that the varieties and stacks (\ref{eq-sec-RGSSM-1}), (\ref{eq-sec-RGSSM-2}) and (\ref{eq-sec-RGSSM-3}) admit integral models; i.e., there exists varieties or stacks that are defined over a finitely generated ring so that base change to $\mathbb{C}$ recovers the original varieties or stacks over $\mathbb{C}$; see \cite[\S 7.1 and \S 7.2.2]{Sco}. A variety or stack over a finite field $\mathbb{F}_q$ is obtained by a similar base change. By abuse of notation, we will denote by $\Rep_B(\mathcal{C})(\mathbb{F}_q)$, $\mathcal{M}_B(\mathcal{C})(\mathbb{F}_q)$ and $\mathfrak{M}_B(\mathcal{C})(\mathbb{F}_q)$ the set or groupoid of $\mathbb{F}_q$-points of the corresponding variety or stack over $\mathbb{F}_q$.

A formula for the number of rational points of generic semi-simple character varieties over $\mathbb{F}_q$ has been obtained by Hausel-Letellier-Rodriguez-Villegas. Let $\mathbf{d}$ be a dimension vector satisfying $\bs\zeta^{\mathbf{d}}=1$, which is supported on the comet-shaped quiver $Q$ as in \S \ref{subsec-multi-comet} and has nonincreasing components along every leg. It defines a $k$-tuples of partitions of $d_{\star}$, each partition consists of the nonzero numbers among $d_{[j,i]}-d_{[j,i+1]}$. We may denoted by $\mathbb{H}_{\mathbf{d}}(\sqrt{q},\frac{1}{\sqrt{q}})$ the rational function in $q$ defined by the $k$-tuple of partitions associated to $\mathbf{d}$.
\begin{Thm}(\cite[Theorem 1.2.3]{HLRV} and \cite[Theorem 5.1.1]{HLRV2})\label{Thm-HLRV}
Suppose that $\mathcal{C}$ is a tuple of generic semi-simple conjugacy classes. Then, we have
$$
|\mathcal{M}_B(\mathcal{C})(\mathbb{F}_q)|=q^{\frac{1}{2}d(\mathcal{M}_{\mathcal{C}})}\mathbb{H}_{\mathbf{d}}(\sqrt{q},\frac{1}{\sqrt{q}}).
$$
Moreover, the polynomial $|\mathcal{M}_B(\mathcal{C})(\mathbb{F}_q)|$ in $q$ is monic.
\end{Thm}

Define
$$
R_{\bs\zeta,\mathbf{d}}:=\{\mathbf{c}\in(\mathbf{Z}_{\ge0})^{Q_0}\mid\bs\zeta^{\mathbf{c}}=1,~\mathbf{c}\le\mathbf{d},\text{ and }c_{[j,i_1]}\ge c_{[j,i_2]}\text{ whenever }i_i\le i_2\}.
$$
\begin{Thm}(\cite[Theorem 1.3.3]{Sco})\label{Thm-Scognamiglio}
For any genus $g$ and any tuple of semi-simple conjugacy classes $\mathcal{C}$, which determine the quiver $Q$, the deformation parameter $\bs\zeta$ and the dimension vector $\mathbf{d}$, we have
$$
|\mathfrak{M}_B(\mathcal{C})(\mathbb{F}_q)|=q^{-(\mathbf{d},\mathbf{d})}\Coeff_{T^{\mathbf{d}}}\big(\Exp(\sum_{\mathbf{c}\in R_{\bs\zeta,\mathbf{d}}}\frac{q\mathbb{H}_{\mathbf{c}}(\sqrt{q},\frac{1}{\sqrt{q}})}{q-1}T^{\mathbf{c}})\big),
$$
where $\Coeff_{T^{\mathbf{d}}}$ means the coefficient of $T^{\mathbf{d}}$ in an element of $\mathbb{Q}(q)[[T_v]]_{v\in Q_0}$.
\end{Thm}

\subsection{Irreducibility via point-counting}\hfill

If the integral model of a complex algebraic variety has polynomial count (that is, for any base change to a finite field $\mathbb{F}_q$, the number of $\mathbb{F}_q$-points is a polynomial in $q$ that is independent of the base change), then the counting polynomial agrees with the E-polynomial of this complex variety. For algebraic stacks, both the E-polynomial and the counting polynomial should be replaced by a rational function; see \cite[Theorem 7.1.2]{Sco}. According to the remarks above the equation \cite[(7.2.6)]{Sco}, character stacks have rational count, and thus counting points gives their E-series.

\begin{Lem}\label{Lem-monicE-irr}
Let $X$ be a complex algebraic variety of pure dimension $d$, and suppose that the coefficient of $x^dy^d$ in its E-polynomial is equal to one. Then, $X$ is irreducible.
\end{Lem}
\begin{proof}
Let $U\subset X$ be the smooth locus of $X$. It suffices to show that $U$ is connected. The E-polynomial of $X\setminus U$ does not contribute to the coefficient of $x^dy^d$. By \cite[Theorem 2.1.3]{HRV}, only the top cohomology of $U$ contains weight $2d$ subquotients, while for any connected component $U_0\subset U$, the one dimensional space $H^{2d}(U_0)$ has pure weight $2d$; therefore, the coefficient of $x^dy^d$ in the E-polynomial of $U$ is equal to the number of the connected components of $U$.
\end{proof}

\begin{Thm}\label{Thm-ani-ss}
Let $Q$, $\bs\zeta$ and $m\mathbf{d}$ be determined by $g$ and $\mathcal{C}$ as in Proposition \ref{Prop-B-open-multi}, so that $\mathcal{M}_B(\mathcal{C})$ is an open subvariety of $\mathcal{M}(\bs\zeta,m\mathbf{d})$, where $m\in\mathbb{Z}_{>0}$ is chosen in such a way that $\mathbf{d}\in R^+_{\bs\zeta}$ is a $\bs\zeta$-indivisible dimension vector. Suppose that $\bs\zeta$ is almost generic for $m\mathbf{d}$ and that $m\mathbf{d}$ lies in $\Sigma_{\bs\zeta}$. Then, the representation variety $\Rep_B(\mathcal{C})$ is irreducible.
\end{Thm}
\begin{proof}
Since $\bs\zeta$ is almost generic for $m\mathbf{d}$, we have $R_{\bs\zeta,m\mathbf{d}}=\{k\mathbf{d}\mid 1\le k\le m\}$. If $m=1$, then $\bs\zeta$ is generic, in which case the irreducibility is known. We will assume $m\ge 2$. Write
$$
f(q,y):=\sum_{1\le k\le m}\frac{q\mathbb{H}_{k\mathbf{d}}(\sqrt{q},\frac{1}{\sqrt{q}})}{q-1}y^{k\mathbf{d}}.
$$
By Theorem \ref{Thm-HLRV}, the Laurent polynomial $q\mathbb{H}_{k\mathbf{d}}(\sqrt{q},\frac{1}{\sqrt{q}})$ is monic of degree $2-k^2(\mathbf{d},\mathbf{d})$, and the lowest degree in it is $k^2(\mathbf{d},\mathbf{d})$. We will write $c_{\mathbf{d}}=-(\mathbf{d},\mathbf{d})$, which satisfies $c_{\mathbf{d}}\ge 1$ since $\mathbf{d}$ is anisotropic (if $\mathbf{d}$ is isotropic, then $m\mathbf{d}\notin\Sigma_{\bs\zeta}$ for $m>1$). We will show that 
\begin{equation}\label{eq-Thm-ani-0}
q^{-m^2(\mathbf{d},\mathbf{d})}|\GL_n(\mathbb{F}_q)|\Coeff_{y^{m\mathbf{d}}}\Exp(f(q,y))
\end{equation}
is a monic polynomial of degree $1+n^2+2m^2c_{\mathbf{d}}$ in $q$. It will then follow from Theorem \ref{Thm-Scognamiglio} that $|\Rep_B(\mathcal{C})(\mathbb{F}_q)|$ is a monic polynomial. 

We have
$$
\Exp(f(q,y))=\sum_{\lambda\in\mathcal{P}}\frac{p_{\lambda}[f(q,y)]}{z_{\lambda}}.
$$
Let $\lambda\in\mathcal{P}$ and write $\lambda=(\lambda_1\ge\lambda_2\ge\cdots\ge\lambda_r)$. Let $\{k_t\ge 1\}_{1\le t\le r}$ be such that $k_1\lambda_1+\cdots+k_r\lambda_r=m$. Recall that
$$
|\GL_n(\mathbb{F}_q)|=q^{\frac{1}{2}n(n-1)}(q-1)(q^2-1)\cdots(q^n-1).
$$ 
Let us show that 
\begin{equation}\label{eq-Thm-ani}
q^{-m^2(\mathbf{d},\mathbf{d})}|\GL_n(\mathbb{F}_q)|\prod_{1\le t\le r}p_{\lambda_t}[\frac{q\mathbb{H}_{k_t\mathbf{d}}(\sqrt{q},\frac{1}{\sqrt{q}})}{q-1}]
\end{equation}
is a polynomial in $q$. Since there is no polar term in $q^{-m^2(\mathbf{d},\mathbf{d})}\prod_{1\le t\le r}p_{\lambda_t}[q\mathbb{H}_{k_t\mathbf{d}}(\sqrt{q},\frac{1}{\sqrt{q}})]$, it suffices to show that $\prod_{1\le t\le r}(q^{\lambda_t}-1)$ is a divisor of $\prod_{1\le t\le n}(q^t-1)$. Define $\mu_1=\lambda_1$. For any $2\le t\le r$, there is a unique integer $\mu_t$ that is divisible by $\lambda_t$ and satisfies $\lambda_1+\cdots+\lambda_{t-1}<\mu_t\le\lambda_1+\cdots+\lambda_{t}$. Now, $\prod_{1\le t\le r}(q^{\lambda_t}-1)$ is a divisor of $\prod_{1\le t\le r}(q^{\mu_t}-1)$, but the latter is a divisor of $\prod_{1\le t\le m}(q^t-1)$ since $\mu_t$'s are distinct and $\mu_r\le m$. The degree of (\ref{eq-Thm-ani}) is equal to $n^2+m^2c_{\mathbf{d}}+\sum_{1\le t\le r}(\lambda_t+\lambda_tk_t^2c_{\mathbf{d}})$; in particular, it has degree $1+n^2+2m^2c_{\mathbf{d}}$ if $\lambda=(1)$ and $k_1=m$. Note that $z_{(1)}=1$. In order to show that (\ref{eq-Thm-ani-0}) is monic of degree $1+n^2+2m^2c_{\mathbf{d}}$, it remains to show that the degree of (\ref{eq-Thm-ani}) is smaller than $1+n^2+2m^2c_{\mathbf{d}}$ either when $r>1$ or when $r=1$ and $\lambda_1>1$. Suppose that $r=1$. Consider
\begingroup
\allowdisplaybreaks
\begin{align}\label{eq-Thm-ani-2}
&1+m^2c_{\mathbf{d}}-\lambda_1-\lambda_1k_1^2c_{\mathbf{d}}\\
\nonumber
=&(c_{\mathbf{d}}\lambda_1k_1^2-1)(\lambda_1-1).
\end{align}
\endgroup 
Note that $c_{\mathbf{d}}\lambda_1k_1^2-1=c_{\mathbf{d}}mk_1-1\ge 1$. The number (\ref{eq-Thm-ani-2}) is zero if $\lambda_1=1$ and $k_1=m$ and is strictly positive otherwise, which proves what we need in the case $r=1$.
Suppose that $r\ge2$ and let us prove by induction on $r$. Consider
\begingroup
\allowdisplaybreaks
\begin{align}\label{eq-Thm-ani-3}
&1+m^2c_{\mathbf{d}}-\sum_{1\le t\le r}(\lambda_t+\lambda_tk_t^2c_{\mathbf{d}})\\
\nonumber
=&1+m_{r-1}^2c_{\mathbf{d}}-\sum_{1\le t\le r-1}(\lambda_t+\lambda_tk_t^2c_{\mathbf{d}})+\lambda_r^2k_r^2c_{\mathbf{d}}+2m_{r-1}\lambda_rk_rc_{\mathbf{d}}-\lambda_r-\lambda_rk_r^2c_{\mathbf{d}}\\
\nonumber
\ge&\lambda_r^2k_r^2c_{\mathbf{d}}+2m_{r-1}\lambda_rk_rc_{\mathbf{d}}-\lambda_r-\lambda_rk_r^2c_{\mathbf{d}}\\
\nonumber
=&\lambda_r\big(c_{\mathbf{d}}k_r(\lambda_rk_r+2m_{r-1}-k_r)-1\big),
\end{align}
\endgroup
where $m_{r-1}=m-k_r\lambda_r$. Now, we have $k_r(\lambda_r-1)+2m_{r-1}\ge2$, and thus (\ref{eq-Thm-ani-3}) is strictly positive. We have thus shown that (\ref{eq-Thm-ani-0}) is monic of degree $1+n^2+2m^2c_{\mathbf{d}}$, and thus $|\Rep_B(\mathcal{C})(\mathbb{F}_q)|$ is a monic polynomial in $q$. 

Finally, we remark that $\Rep_B(\mathcal{C})$ is pure dimensional. It then follows from Lemma \ref{Lem-monicE-irr} that $\Rep_B(\mathcal{C})$ is irreducible. Since $m\mathbf{d}\in\Sigma_{\bs\zeta}$, the representation variety of quiver representations $\Rep^{\circ}(\bar{Q},\bs\zeta,m\mathbf{d})$, where $Q$, $\bs\zeta$ and $\mathbf{d}$ are as in Proposition \ref{Prop-B-open-multi}, is pure dimensional according to \cite[Theorem 1.11]{CBS}. In the proof of Proposition \ref{Prop-B-open-multi}, we have identified $\Rep_B(\mathcal{C})$ with an open subset of a variety $\Rep^{\circ}(g,\mathcal{C})$, which turns out to be isomorphic to the quotient
$$
\Rep^{\circ}(\bar{Q},\bs\zeta,m\mathbf{d})\ds G_{>0};
$$
see (\ref{eq-Prop-B-open-multi-4}). Since we assume $C_j$'s to be semi-simple conjugacy classes, we have $\mathcal{C}=\bar{\mathcal{C}}$. Therefore, for every $(\mathbf{V},\bs\phi)\in\Rep^{\circ}(\bar{Q},\bs\zeta,m\mathbf{d})$, the components $\phi^{\ast}_{[j,i]}$ are all injective. This implies that $G_{>0}$ acts freely on $\Rep^{\circ}(\bar{Q},\bs\zeta,m\mathbf{d})$, and thus the quotient variety is pure dimensional.
\end{proof}

\section{The proof of the main theorem}\label{sec-Main}

Our ultimate goal is to prove that a general Betti moduli space $\mathcal{M}_B(\mathbf{P},\bar{\mathcal{C}},\bs\beta)$, where neither $\bs\beta$ nor $\bar{\mathcal{C}}$ needs to be generic, is connected. The core of our proof consists of producing a transcendental surjective map from a moduli space $\mathcal{M}_B(\mathcal{C}')$, which is known to be connected, to $\mathcal{M}_B(\mathbf{P},\bar{\mathcal{C}},\bs\beta)$. To this end, we use the nonabelian Hodge correspondence to pass to the moduli spaces of parabolic Higgs bundles. We will give a brief review of Higgs bundles as well as the statement of the nonabelian Hodge correspondence in \S \ref{subsec-NAH}. The key fact is that a modification of the parabolic weight (the stability condition) on the Dolbeault side corresponds to a modification of the eigenvalues of the monodromy (the deformation parameter) on the Betti side. More details of the desired transcendental map will be given in \S \ref{subsec-diagram}, and the proof will be completed in \S \ref{subsec-main-proof}.

\subsection{Nonabelian Hodge correspondence}\label{subsec-NAH}\hfill

We regard $\bar{X}$ as a smooth algebraic curve and $\{x_1,\ldots,x_k\}$ as a divisor $D$. A meromorphic Higgs bundle on $\bar{X}$ is a pair $(E,\Phi)$, where $E$ is a locally free sheaf on $\bar{X}$ and $\Phi:E\rightarrow E\otimes\Omega(D)$ is a homomorphism of coherent sheaves, where $\Omega$ is the canonical bundle of $\bar{X}$. We will denote by $\Phi_j$ the residue of the Higgs field at $x_j$, which is regarded as an endomorphism of the fibre $E_{x_j}$. A parabolic structure on the Higgs bundle $(E,\Phi)$ at $x_j$ is the data of a flag $E_{[j,\bullet]}$ on $E_{x_j}$ that is preserved by $\Phi_j$. We will represent a parabolic structure by a parabolic Lie subalgebra $\mathfrak{p}_j\subset\mathfrak{gl}(E_{x_j})$, and denote by $\mathfrak{l}_j$ its Levi quotient. A parabolic Higgs bundle is a meromorphic Higgs bundle equipped with a parabolic structure; we will denote by $(E,E_{[j,\bullet]},\Phi)$ a parabolic Higgs bundle. Define a tuple of integers $(d_{[j,s]})_{1\le j\le k,~s\ge0}$ by $d_{[j,s]}=\dim E_{[j,s]}$. For each $j$, there is an integer $\mu_j$ such that $d_{[j,\mu_j]}>0$ and $d_{[j,\mu_j+1]}=0$. For a given parabolic structure on $E$, we will denote by $\bar{\Phi}_j$ the induced element of 
$$
\mathfrak{l}_j\cong\bigoplus_{s=0}^{\mu_j}\gl(E_{[j,s]}/E_{[j,s+1]}),
$$ 
which will be called the graded residue of $\Phi$ at $x_j$. A parabolic weight for $(E,E_{[j,\bullet]},\Phi)$ is a tuple of real numbers $\bs\alpha=(\bs\alpha_{[j,s]})_{1\le j\le k,~0\le s\le \mu_j}$ satisfying
$$
0\le \alpha_{[j,0]}<\alpha_{[j,1]}<\cdots<\alpha_{[j,\mu_j]}<1.
$$
The parabolic degree of $\mathscr{E}=(E,E_{[j,\bullet]},\Phi)$ with weight $\bs\alpha$ is 
$$
\deg_{\bs\alpha}\mathscr{E}:=\deg E+\sum_{j=1}^k\sum_{s=0}^{\mu_j}\alpha_{[j,s]}(d_{[j,s]}-d_{[j,s+1]}).
$$
We say that $(E,E_{[j,\bullet]},\Phi)$ is $\bs\alpha$-semi-stable (resp. $\bs\alpha$-stable) if for any vector subbundle $F\subset E$ that is preserved under $\Phi$, the induced parabolic Higgs subbundle $\mathscr{F}=(F,F_{[\bullet,\bullet]},\Phi)$ with the induced parabolic weight $\bs\alpha$ satisfies
$$
\frac{\deg_{\bs\alpha}\mathscr{F}}{\rk F}\le \frac{\deg_{\bs\alpha}\mathscr{E}}{\rk E}\quad\text{resp. }\frac{\deg_{\bs\alpha}\mathscr{F}}{\rk F}< \frac{\deg_{\bs\alpha}\mathscr{E}}{\rk E}.
$$
Write $\mathbf{p}=(\mathfrak{p}_j)_{1\le j\le k}$. If the parabolic structure of a parabolic Higgs bundle $(E,E_{[j,\bullet]},\Phi)$ at $x_j$ is conjugate to a given parabolic subalgebra $\mathfrak{p}_j$ for every $j$, then we will say that it is of parabolic type $\mathbf{p}$. Let $\mathcal{O}=(\mathcal{O}_j)_{1\le j\le k}$ be a semi-simple adjoint orbit of $\prod_{j=1}^k\mathfrak{l}_j$, and let $\delta\in\mathbb{Z}$. There exists a coarse moduli space 
$$
\mathcal{M}_{Dol}(\mathbf{p},\mathcal{O},\bs\alpha,\delta)
$$
parametrising $\bs\alpha$-semi-stable parabolic Higgs bundles of degree $\delta$ and parabolic degree zero and of parabolic type $\mathbf{p}$ whose Higgs field at $x_j$ has a graded residue lying in $O_j$ for each $j$.

The nonabelian Hodge correspondence (see \cite[Theorem 7.10]{BGM} or \cite[Theorem 5.3]{HS}) implies that there is a homeomorphism of topological spaces
\begin{equation}\label{eq-NAH}
\mathcal{M}_B(\mathbf{P},\mathcal{C},\bs\beta)\longrightarrow\mathcal{M}_{Dol}(\mathbf{p},\mathcal{O},\bs\alpha,\delta),
\end{equation}
where $\mathcal{C}$ is a tuple of semi-simple conjugacy classes, and the invariants on the two sides (i.e. $\mathbf{P}$, $\mathcal{C}$ and $\bs\beta$ on the Betti side, and $\mathbf{p}$, $\mathcal{O}$ and $\bs\alpha$ on the Dolbeault side) determine each other; see \cite[\S 5]{Si90}, \cite[Table 1]{BGM}, \cite[\S 5]{HS} or \cite[\S 3.4]{Shu25} for the precise formulae. Note that the map (\ref{eq-NAH}) has tacitly incorporated the Riemann-Hilbert correspondence.

Similar to \S \ref{subsec-KM-roots}, we have the notion of almost generic parabolic weights. Fix a tuple of integers $(\delta,(d_{[j,s]})_{1\le j\le k,~0\le s\le \mu_j})$, where $d_{[j,0]}=n$ for all $j$ and $d_{[j,s]}>d_{[j,s+1]}$ for all $j$ and $s$. Let $\bs\alpha=(\bs\alpha_{[j,s]})_{1\le j\le k,~0\le s\le \mu_j}$ be a tuple of real numbers satisfying
$$
\delta+\sum_{j=1}^k\sum_{s=0}^{\mu_j}\alpha_{[j,s]}(d_{[j,s]}-d_{[j,s+1]})=0.
$$ 
We say that $\bs\alpha$ is almost generic if every tuple integers $(\delta',(d'_{[j,s]})_{1\le j\le k,~0\le s\le \mu_j})$ that satisfies
\begin{itemize}
\item $d'_{[j,0]}=n'$ for some $n'\le n$ and all $j$,
\item $d'_{[j,s]}>d'_{[j,s+1]}$ and $d'_{[j,s]}\le d_{[j,s]}$ for all $j$ and $s$, and
\item the equality 
$$
\delta'+\sum_{j=1}^k\sum_{s=0}^{\mu_j}\alpha_{[j,s]}(d'_{[j,s]}-d'_{[j,s+1]})=0,
$$ 
\end{itemize} 
is a $\mathbb{Q}$-multiple of $(\delta,(d_{[j,s]})_{1\le j\le k,~0\le s\le \mu_j})$; if such a tuple $(\delta',(d'_{[j,s]})_{1\le j\le k,~0\le s\le \mu_j})$ is either zero or $(\delta,(d_{[j,s]})_{1\le j\le k,~0\le s\le \mu_j})$ itself, then we say that $\bs\alpha$ is generic.

\subsection{A diagram}\label{subsec-diagram}\hfill

We will use the following diagram to produce a continuous surjective map from $\mathcal{M}_B(\mathcal{C}')$ to $\mathcal{M}_B(\mathbf{P},\bar{\mathcal{C}},\bs\beta)$; it is based on a construction of Ballandras \cite{Ball} and was used in \cite{Shu25} for the same purpose.
\begin{equation}\label{eq-red-diag}
\begin{tikzcd}[row sep=2.5em, column sep=2.5em]
\mathcal{M}_B(\mathbf{P},\bar{\mathcal{C}},\bs\beta) & \mathcal{M}_B(\mathbf{P}^{\sharp},\mathcal{C}^{\sharp},\bs\beta') \arrow[l, "\mathbf{Var}_1"'] \arrow[r, "\mathbf{NH}_1"] & \mathcal{M}_{Dol}(\mathbf{p}',\mathcal{O}',\bs\alpha') \\
\mathcal{M}_B(\mathcal{C}') & \mathcal{M}_B(\mathbf{P}',\mathcal{C}',\bs\beta') \arrow[l, "\mathbf{Iso}"] \arrow[r, "\mathbf{NH}_2"'] & \mathcal{M}_{Dol}(\tilde{\mathbf{p}},\mathcal{O},\bs\alpha). \arrow[u, "\mathbf{Var}_2"']
\end{tikzcd}
\end{equation}
Our object of interest is $\mathcal{M}_B(\mathbf{P},\bar{\mathcal{C}},\bs\beta)$, which we would like to prove to be connected. The connectedness of $\mathcal{M}_B(\mathcal{C}')$, as we will explain, is what we proved in \S \ref{sec-RGSSM}. Ideally, the maps $\mathbf{Var}_1$ and $\mathbf{Var}_2$ are surjective and the maps $\mathbf{NH}_1$, $\mathbf{NH}_2$ and $\mathbf{Iso}$ are bijective, and thus every space in this diagram is connected. However, the surjectivity of $\mathbf{Var}_1$ and $\mathbf{Var}_2$ is conditional on the type of $\mathbf{P}$ or $\mathbf{p}'$. We will leave the clarification of the surjectivity problem to the next subsection; for now, we explain in more detail the constituents of this diagram.

The maps $\mathbf{NH}_1$ and $\mathbf{NH}_2$ are nonabelian Hodge correspondences, and only the invariants defining the relevant moduli spaces need to be clarified. We start from the space $\mathcal{M}_B(\mathbf{P},\bar{\mathcal{C}},\bs\beta)$, and regard it as an open subvariety of a multiplicative quiver variety $\mathcal{M}(\bs\zeta,\mathbf{d},\bs\theta)$. We choose $\bs\theta'$ to be an almost generic stability condition in a small neighbourhood of $\bs\theta$ satisfying $\theta'_{[j,i]}>0$ for any $j$ and $i>0$. Then, variation of stability conditions defines a projective morphism
\begin{equation}\label{eq-Var-theta'}
\mathbf{Var}_1:\mathcal{M}(\bs\zeta,\mathbf{d},\bs\theta')\longrightarrow\mathcal{M}(\bs\zeta,\mathbf{d},\bs\theta).
\end{equation}
We claim that the inverse image of $\mathcal{M}_B(\mathbf{P},\bar{\mathcal{C}},\bs\beta)$ is a Betti moduli space. In the proof of Proposition \ref{Prop-B-open-multi}, we have shown that there is an isomorphism
$$
\Rep^{\circ}(\bar{Q},\bs\zeta,\mathbf{d})^{\bs\theta-ss}\ds G_{>0}\cong\Rep^{\circ}(g,\mathbf{P},\bar{\mathcal{C}})^{\bs\beta-ss},
$$
and that there is an open subset $U\subset\Rep^{\circ}(g,\mathbf{P},\bar{\mathcal{C}})$ such that 
$$
\mathcal{M}_B(\mathbf{P},\bar{\mathcal{C}},\bs\beta)\cong(U\cap \Rep^{\circ}(g,\mathbf{P},\bar{\mathcal{C}})^{\bs\beta-ss})\ds\GL_n.
$$
By the choice of $\bs\theta'$, we have an inclusion $\Rep^{\circ}(\bar{Q},\bs\zeta,\mathbf{d})^{\bs\theta'-ss}\subset\Rep^{\circ}(\bar{Q},\bs\zeta,\mathbf{d})^{\bs\theta-ss}$, which induces
$$
\Rep^{\circ}(\bar{Q},\bs\zeta,\mathbf{d})^{\bs\theta'-ss}\ds G_{>0}\longrightarrow\Rep^{\circ}(\bar{Q},\bs\zeta,\mathbf{d})^{\bs\theta-ss}\ds G_{>0},
$$
and the quotient of this map by $\GL_n$ is $\mathbf{Var}_1$. Now, the same arguments as in the proof of Proposition \ref{Prop-B-open-multi} give an isomorphism
$$
\Rep^{\circ}(\bar{Q},\bs\zeta,\mathbf{d})^{\bs\theta'-ss}\ds G_{>0}\cong\Rep^{\circ}(g,\mathbf{P}^{\sharp},\mathcal{C}^{\sharp})^{\bs\beta'-ss},
$$
for a tuple of parabolic subgroups $\mathbf{P}^{\sharp}$, a tuple of conjugacy classes $\mathcal{C}^{\sharp}$ and a tuple of  parabolic weights $\bs\beta'$. Note that every component of $\mathbf{P}^{\sharp}=(P_j^{\sharp})_{1\le j\le k}$ is a parabolic subgroup that is contained in the corresponding component of $\mathbf{P}$ and that $\mathcal{C}^{\sharp}$ consists of central conjugacy classes of the Levi quotients of the $P_j^{\sharp}$'s (this is precisely the case studied by Yamakawa \cite{Yama}). It follows from the definition of $U$ (see (\ref{eq-Prop-B-open-multi-U})) that its inverse image in $\Rep^{\circ}(g,\mathbf{P}^{\sharp},\mathcal{C}^{\sharp})^{\bs\beta'-ss}$ is the open subset $U'$ defined by the invertibility of the maps $(\psi_i)_{1\le i\le g}$. We thus see that $(U'\cap \Rep^{\circ}(g,\mathbf{P}^{\sharp},\mathcal{C}^{\sharp})^{\bs\beta'-ss})\ds\GL_n$ is the inverse image of $\mathcal{M}_B(\mathbf{P},\bar{\mathcal{C}},\bs\beta)$ and is a Betti moduli space defined by $g$, $\mathbf{P}^{\sharp}$, $\mathcal{C}^{\sharp}$ and $\bs\beta'$. This shows that $\mathbf{Var}_1$ restricts to
$$
\mathbf{Var}_1:\mathcal{M}_B(\mathbf{P}^{\sharp},\mathcal{C}^{\sharp},\bs\beta')\longrightarrow\mathcal{M}_B(\mathbf{P},\bar{\mathcal{C}},\bs\beta).
$$ 
The map $\mathbf{Iso}$ is defined similarly, which turns out to be an isomorphism as we will explain later.

Variation of parabolic weights defines the map $\mathbf{Var}_2$; however, we need to simultaneously alter other invariants as follows. Let $(E,E'_{[\bullet,\bullet]},\Phi)\in\mathcal{M}_{Dol}(\mathbf{p}',\mathcal{O}',\bs\alpha')$. The parabolic structure at $x_j$ is a flag
$$
E_{x_j}=E'_{[j,0]}\supset E'_{[j,1]}\supset\cdots\supset E'_{[j,\mu'_j]}\supset E'_{[j,\mu'_j+1]}=0,
$$
and the endomorphisms of $E_{x_j}$ preserving this flag form the parabolic Lie subalgebra $\mathfrak{p}'_j$, the $j$-th component of $\mathbf{p}'$. Each $\mathcal{O}'_j$ is a semi-simple adjoint orbit in $\bigoplus_{s=0}^{\mu'_j}\gl(E'_{[j,s]}/E'_{[j,s+1]})$, and the graded residue $\bar{\Phi}_j$ lies in $\mathcal{O}'_j$. For any $[j,s]$, we choose an ordered tuple $(\lambda_0,\lambda_1,\ldots,\lambda_{w_s})$ consisting of the distinct eigenvalues of the $E'_{[j,s]}/E'_{[j,s+1]}$-component of $\bar{\Phi}_j$. Then, there is a unique $\bar{\Phi}_j$-invariant flag
$$
E'_{[j,s]}=E_{[j,s],0}\supset E_{[j,s],1}\supset\cdots\supset E_{[j,s],w_s+1}=E'_{[j,s+1]}
$$ 
such that $\bar{\Phi}_j$ acts as the scalar $\lambda_t$ on $E_{[j,i],t}/E_{[j,i],t+1}$. In this way, we find a flag $E_{[\bullet,\bullet]}$ refining $E'_{[\bullet,\bullet]}$, and the endomorphisms of $E_{x_j}$ preserving this flag form a parabolic Lie subalgebra $\tilde{\mathfrak{p}}_j\subset\mathfrak{p}_j'$; this defines $\tilde{\mathbf{p}}$. The residue $\Phi_j$ acts as a scalar on every component of  the associated graded of $E_{[\bullet,\bullet]}$, and thus determines a central adjoint orbit $\mathcal{O}_j$ in the Levi quotient of $\tilde{\mathfrak{p}}_j$. The parabolic weight $\bs\alpha$ is chosen to be almost generic and to lie in a small neighbouhood of $\bs\alpha'$. It follows from \cite[Lemma 3.11 and Proposition 3.12]{Shu25} that $\mathbf{Iso}$ is indeed an isomorphism.
\begin{Rem}\label{Rem-unchange-d}
It can be seen from the transformation rule \cite[Equation (3.4.1)]{Shu25} that if we apply Proposition \ref{Prop-B-open-multi} to $\mathcal{M}_B(\mathbf{P}',\mathcal{C}',\bs\beta')$, then the resulting dimension vector $\mathbf{d}'$ is in fact equal to $\mathbf{d}$. The reason is: the parabolic weights $\beta'_{[j,i]}$ on the Betti side of $\mathbf{NH}_1$ become part of the eigenvalues on the Dolbeault side, and thus are not affected by $\mathbf{Var}_2$; when we return to the Betti side via $\mathbf{NH}_2$, they are still the parabolic weights as in $\mathcal{M}_B(\mathbf{P}^{\sharp},\mathcal{C}^{\sharp},\bs\beta')$ with the same order and the same multiplicities. Moreover, we have in fact $\mathbf{P}^{\sharp}=\mathbf{P}'$.
\end{Rem}

The main property we need about $\mathbf{Var}_2$ is the following.
\begin{Prop}\label{Prop-Var2-birat}
The morphism $\mathbf{Var}_2$ is projective and for any open subset $U\subset\mathcal{M}_{Dol}(\mathbf{p}',\mathcal{O}',\bs\alpha')$   consisting of $\bs\alpha'$-stable parabolic Higgs bundles, the restriction of $\mathbf{Var}_2$ to $\mathbf{Var}_2^{-1}(U)$ is an isomorphism. 
\end{Prop}
\begin{proof}
The projectivity comes from the fact that $\mathbf{Var}_2$ is defined by a variation of stability conditions. The second statement follows from \cite[Proposition 3.8]{Shu25} and the fact that stable loci are smooth.
\end{proof}

\subsection{The proof of the main theorem}\label{subsec-main-proof}\hfill

It is crucial to use the multiplicative quiver variety point of view in our discussion of the connectedness of character varieties. Of particular importance is the type of the dimension vector.

\begin{Prop}\label{Prop-red-imaroots}
Suppose that $\mathcal{M}_B(\mathbf{P},\bar{\mathcal{C}},\bs\beta)$ is nonempty. Then, there is an isomorphism
$$
\mathcal{M}_B(\mathbf{P},\bar{\mathcal{C}},\bs\beta)\lisom \mathcal{M}_B(\mathbf{P}^{\flat},\bar{\mathcal{C}}^{\flat},\bs\beta^{\flat}),
$$
for some $\mathbf{P}^{\flat}$, $\bar{\mathcal{C}}^{\flat}$, $\bs\beta^{\flat}$, such that the quiver $Q^{\flat}$ and the dimension vector $\mathbf{d}^{\flat}$ determined by these data are such that $\mathbf{d}^{\flat}$ is an imaginary root.
\end{Prop}
\begin{proof}
If $g=0$, so that $\mathcal{M}_B(\mathbf{P},\bar{\mathcal{C}},\bs\beta)=\mathcal{M}(Q,\bs\zeta,\mathbf{d},\bs\theta)$, then we apply the standard quiver variety techniques as in the proof of \cite[Theorem 4.6]{Shu25}. For the cases $g>0$, we will follow the same line of arguments and at the same time show that various operations on $\mathcal{M}(Q,\bs\zeta,\mathbf{d},\bs\theta)$ preserve the Betti moduli space contained in it. 

Suppose that $v$ is a vertex with $(\mathbf{d},\mathbf{e}_v)>0$ and either $\zeta_v\ne1$ or $\theta_v\ne0$. Note that there is no loop at $v$, since $(\mathbf{d},\mathbf{e}_v)>0$. Write $\bs\zeta^{\ast}=u_v(\bs\zeta)$, $\mathbf{d}^{\ast}=s_v(\mathbf{d})$ and $\bs\theta^{\ast}=r_v(\bs\theta)$. The reflection functor at $v$ induces an isomorphism
$$
\mathcal{M}(Q,\bs\zeta,\mathbf{d},\bs\theta)\lisom\mathcal{M}(Q,\bs\zeta^{\ast},\mathbf{d}^{\ast},\bs\theta^{\ast}).
$$
If $\theta_v>0$, then $\theta^{\ast}_v<0$; therefore, we cannot immediately identify an open subset of $\mathcal{M}(Q,\bs\zeta^{\ast},\mathbf{d}^{\ast},\bs\theta^{\ast})$ with a Betti moduli space. We can nevertheless define an open subvariety
$$
\mathcal{M}_B(Q,\bs\zeta^{\ast},\mathbf{d}^{\ast},\bs\theta^{\ast})\subset\mathcal{M}(Q,\bs\zeta^{\ast},\mathbf{d}^{\ast},\bs\theta^{\ast})
$$ 
by the invertibility condition on the arrows associated to the loops as in (\ref{eq-Prop-B-open-multi-U}), which will also be called Betti moduli spaces in what follows. This open subset matches $\mathcal{M}_B(Q,\bs\zeta,\mathbf{d},\bs\theta)$ under the above isomorphism. Indeed, the reflection isomorphism by definition does not alter the components of a quiver representation corresponding to the loops. 

If $v$ is a vertex with $(\mathbf{d},\mathbf{e}_v)>0$, $\zeta_v=1$ and $\theta_v=0$, then any $\mathbf{d}$-dimensional $\bs\theta$-polystable representation has a direct factor with dimension vector $\mathbf{e}_v$; this is simply the nontrivial stability version of \cite[Lemma 5.1]{CBS}, the proof of which remains valid in this generality. It follows that taking the direct sum with the unique simple representation of dimension vector $\mathbf{e}_v$ gives an isomorphism
$$
\mathcal{M}(Q,\bs\zeta,\mathbf{d}-\mathbf{e}_v,\bs\theta)\lisom\mathcal{M}(Q,\bs\zeta,\mathbf{d},\bs\theta).
$$
Since the Betti moduli spaces are open subsets defined by invertibility conditions on the maps associated to loops, and there is no loop at $v$, the above isomorphism restricts to Betti moduli spaces.

By induction on $|\mathbf{d}|:=\sum_{v\in Q_0}d_v$, we arrive at a variety $\mathcal{M}(Q,\bs\zeta^{\natural},\mathbf{d}^{\natural},\bs\theta^{\natural})$ with $\mathbf{d}^{\natural}$ satisfying $(\mathbf{d}^{\natural},\mathbf{e}_v)\le0$ for every $v\in Q_0$. Let $Q^{\natural}\subset Q$ be the support of $\mathbf{d}^{\natural}$. The restriction of $\mathbf{d}^{\natural}$ to each connected component of $Q^{\natural}$ lies in the fundamental region. Since $Q$ is comet-shaped and there exists no imaginary root supported on a type $A$ quiver, $Q^{\natural}$ is connected and contains the vertex $\star$. Now, it is possible that $\bs\theta^{\natural}$ is negative at some vertices $[j,i]\ne\star$, so that $\mathcal{M}_B(Q^{\natural},\bs\zeta^{\natural},\mathbf{d}^{\natural},\bs\theta^{\natural})$ is not visibly isomorphic to a tame character variety. Fortunately, the algorithm at the beginning of the proof of \cite[Theorem 4.1]{Shu25} gives a sequence of reflections that take $\mathcal{M}(Q^{\natural},\bs\zeta^{\natural},\mathbf{d}^{\natural},\bs\theta^{\natural})$ to some isomorphic $\mathcal{M}(Q^{\natural},\bs\zeta'',\mathbf{d}'',\bs\theta'')$, with $\bs\theta''$ being nonnegative on every vertex $v\ne \star$. Since $\mathbf{d}^{\natural}$ lies in the fundamental region, the vector $\mathbf{d}''$ obtained by reflections is thus an imaginary root. The desired $\mathcal{M}_B(\mathbf{P}^{\flat},\bar{\mathcal{C}}^{\flat},\bs\beta^{\flat})$ is the Betti moduli space contained in $\mathcal{M}(Q^{\natural},\bs\zeta'',\mathbf{d}'',\bs\theta'')$.
\end{proof}

Now, we discuss the surjectivity issue in the diagram (\ref{eq-red-diag}). The proof of the above proposition says that there is a sequence of reflections taking $\mathbf{d}''$ to a vector $\mathbf{d}^{\natural}$ in the fundamental region. According to \cite[Theorem 6.16]{ST}, three possibilities occur:
\begin{itemize}
\item[($\mathbf{Aff}$)] $\mathbf{d}^{\natural}=m\mathbf{d}_0$, where $m\in\mathbb{Z}_{>1}$ and $\mathbf{d}_0$ is a $\bs\zeta$-indivisible isotropic root supported on an affine Dynkin quiver $\tilde{Q}\subset Q$.
\item[($\mathbf{Aff}^{\infty}$)] $\mathbf{d}^{\natural}=\mathbf{e}_{\infty}+m\mathbf{d}_0$, where $m\in\mathbb{Z}_{>1}$, $\mathbf{d}_0$ is an indivisible isotropic root supported on an affine Dynkin quiver $\tilde{Q}$, $\mathbf{e}_{\infty}$ is the coordinate vector of an additional vertex $\infty$ joined to the affine node of $\tilde{Q}$ by a single arrow.
\item $\mathbf{d}^{\natural}\in\Sigma_{\bs\zeta^{\natural},\bs\theta^{\natural}}$
\end{itemize}
Note that there is another case listed in the aforementioned theorem of Schedler-Tirelli; namely, the supporting quiver of $\mathbf{d}$ consists of two disjoint subsets of vertices $K\sqcup J\subset Q_0$, which are joined by a single edge connecting $v_1\in K$ and $v_2\in J$, while $d_{v_1}=d_{v_2}=1$. However, this case is trivial for star-shaped quiver, since we can always make the components of $\mathbf{d}$ strictly decreasing along every leg of $Q$ by performing a sequence of reflections (an additional real root supported on a leg does not affect the isomorphism type of the resulting quiver variety, and thus can be omitted). 

Indeed, the surjectivity of $\mathbf{Var}_1$ and $\mathbf{Var}_2$ can be proved for the roots in $\Sigma_{\bs\zeta,\bs\theta}$. For the other two types of dimension vectors, the connectedness of character varieties can be proved by a different method. The following proposition allows $g>0$ and thus is a generalisation of \cite[Proposition 4.2]{Shu25}.

\begin{Prop}\label{Prop-conn-zeta-indiv}
Let $\mathcal{M}_B(\mathbf{P},\bar{\mathcal{C}},\bs\beta)$ be a nonempty Betti moduli space and let $\mathcal{M}(\bs\zeta,\mathbf{d},\bs\theta)$ be the multiplicative quiver variety as in Proposition \ref{Prop-B-open-multi}, which contains $\mathcal{M}_B(\mathbf{P},\bar{\mathcal{C}},\bs\beta)$ as an open subset. Suppose that a sequence of reflections take $\mathbf{d}$ to $\mathbf{d}^{\natural}$, which lies in the fundamental region, an that $\mathbf{d}^{\natural}$ is of type $(\mathbf{Aff})$ or $(\mathbf{Aff}^{\infty})$. Then, the variety $\mathcal{M}_B(\mathbf{P},\bar{\mathcal{C}},\bs\beta)$ is connected.
\end{Prop}
\begin{proof}
If $g=0$, then the proposition follows from \cite[Proposition 4.2 and Theorem 4.4]{Shu25}. Suppose $g=1$. Then, the case $(\mathbf{Aff})$ is precisely the character varieties for elliptic curves with a central monodromy, which are known to be connected; see for example \cite[\S 6.2 and \S 6.3]{LNY} or \cite[\S 7.1]{Shu24}. For the case $(\mathbf{Aff}^{\infty})$, note that the same proof for $g=0$ remains valid. Finally, there is no dimension vector of either type if $g\ge2$.
\end{proof}

\begin{Thm}\label{Thm-Irr-central-mono}
The moduli space $\mathcal{M}_B(\mathbf{P},\bar{\mathcal{C}},\bs\beta)$, if nonempty, is an irreducible algebraic variety.
\end{Thm}
\begin{proof}
Let $Q$, $\bs\zeta$, $\mathbf{d}$ and $\bs\theta$ be the quiver data determined by $g$, $\mathbf{P}$, $\bar{\mathcal{C}}$, $\bs\beta$. According to the proof of Proposition \ref{Prop-red-imaroots}, we may assume that $\mathbf{d}$ is an imaginary root (or rather, we may replace $(Q,\mathbf{d})$ by $(Q^{\natural},\mathbf{d}'')$), and that there is a sequence of reflections taking $\mathcal{M}(Q,\bs\zeta,\mathbf{d},\bs\theta)$ to $\mathcal{M}(Q^{\natural},\bs\zeta^{\natural},\mathbf{d}^{\natural},\bs\theta^{\natural})$, where $\mathbf{d}^{\natural}$ lies in the fundamental region. It suffices to prove the irreducibility in the following cases: (1) $\mathbf{d}^{\natural}$ is of type $(\mathbf{Aff})$. (2) $\mathbf{d}^{\natural}$ is of type $(\mathbf{Aff}^{\infty})$. (3) $\mathbf{d}^{\natural}\in\Sigma_{\bs\zeta^{\natural},\bs\theta^{\natural}}$. Cases (1) and (2) are precisely the content of Proposition \ref{Prop-conn-zeta-indiv}. Now we assume $\mathbf{d}^{\natural}$ to be of the form (3). Since $\mathbf{d}$ is obtained from $\mathbf{d}^{\natural}$ by reflections, we have $\mathbf{d}\in\Sigma_{\bs\zeta,\bs\theta}$. 

Consider the diagram (\ref{eq-red-diag}). According to Remark \ref{Rem-unchange-d}, the dimension vector associated to $\mathcal{M}_B(\mathbf{P}',\mathcal{C}',\bs\beta')$ is still $\mathbf{d}$ supported on the same quiver $Q$. Let $\bs\zeta'$, $\mathbf{d}$ and $\bs\theta'$ be the quiver data determined by $\mathbf{P}'$, $\mathcal{C}'$ and $\bs\beta'$. If we write $\mathbf{d}=\sum_{s=1}^r\mathbf{d}^{(s)}$ with $\mathbf{d}^{(s)}\in R^+_{\bs\zeta',\bs\theta'}$, then each $\mathbf{d}^{(s)}$ is a $\mathbb{Q}$-multiple of $\mathbf{d}'$. But $\mathbf{d}^{(s)}$ also lies in $R^+_{\bs\zeta,\bs\theta'}$ by the definition of $\bs\zeta'$. It follows that $\mathbf{d}\in\Sigma_{\bs\zeta',\bs\theta'}=\Sigma_{\bs\zeta'}$. and thus $\mathcal{M}_B(\mathcal{C}')$ satisfies  the assumptions of Theorem \ref{Thm-ani-ss}.

We will show that the diagram (\ref{eq-red-diag}) gives a surjective morphism $\mathcal{M}_B(\mathcal{C}')\rightarrow \mathcal{M}_B(\mathbf{P},\bar{\mathcal{C}},\bs\beta)$, from which we deduce that $\mathcal{M}_B(\mathbf{P},\bar{\mathcal{C}},\bs\beta)$ is connected. Since $\mathcal{M}_B(\mathbf{P},\bar{\mathcal{C}},\bs\beta)$ is also known to be normal by \cite[Theorem 5.4]{ST} and \cite{CB03}, the irreducibility follows. To prove surjectivity, first note that the map $\mathbf{Iso}$ is an isomorphism, as we have explained. By assumption, we have $\mathbf{d}\in\Sigma_{\bs\zeta,\bs\theta}$, and thus $\mathcal{M}_B(\mathbf{P},\bar{\mathcal{C}},\bs\beta)^s\subset \mathcal{M}_B(\mathbf{P},\bar{\mathcal{C}},\bs\beta)$ is an open dense subset by Corollary \ref{Cor-dense-M_B^s}. Since $\mathbf{Var}_1$ is projective, it has to be surjective. It then remains to show that the image of $\mathbf{Var}_2$ contains an open dense subset, noting that a variation of parabolic weights induces a projective morphism. We claim that $\Sigma_{\bs\zeta,\bs\theta}\subset\Sigma_{\bs\zeta,\bs\theta'}$ by the definition of these sets, where $\bs\theta'$ is as in (\ref{eq-Var-theta'}). Indeed, we have $R^+_{\bs\zeta,\bs\theta'}\subset R^+_{\bs\zeta,\bs\theta}$; this is because $\bs\theta'$ is almost generic, and thus elements of $R^+_{\bs\zeta,\bs\theta'}$ are $\mathbb{Q}$-multiples of $\mathbf{d}$, which necessarily lie in $R^+_{\bs\zeta,\bs\theta}$. It follows that $\mathcal{M}_B(\mathbf{P}^{\sharp},\mathcal{C}^{\sharp},\bs\beta')^s\subset \mathcal{M}_B(\mathbf{P}^{\sharp},\mathcal{C}^{\sharp},\bs\beta')$ is an open dense subset. Since $\mathbf{NH}_1$ is a homeomorphism, the loci of stable parabolic Higgs bundles is dense in $\mathcal{M}_{Dol}(\mathbf{p}',\mathcal{O}',\bs\alpha')$. Note that the image of $\mathbf{Var}_2$ contains this subset by Proposition \ref{Prop-Var2-birat}. We have thus proved the surjectivity of $\mathbf{Var}_2$, hence the irreducibility of $\mathcal{M}_B(\mathbf{P},\bar{\mathcal{C}},\bs\beta)$.
\end{proof}

\section{Symplectic leaves}\label{sec-Symp-L}

Our study of the symplectic leaves of character varieties is based on the quiver variety point of view, and we will always embed $\mathcal{M}_B(\mathbf{P},\bar{\mathcal{C}},\bs\beta)$ as an open subvariety in some $\mathcal{M}(\bs\zeta,\mathbf{d},\bs\theta)$, and write
$$
\mathcal{M}_B(\bs\zeta,\mathbf{d},\bs\theta):=\mathcal{M}_B(\mathbf{P},\bar{\mathcal{C}},\bs\beta),
$$ 
tacitly using Proposition \ref{Prop-B-open-multi}. We will classify the symplectic leaves and the codimension-two leaves in tame character varieties in Theorem \ref{Thm-rep-str=leaves} and Proposition \ref{Prop-codim2-leaves}, which are immediate consequences of our connectedness theorem and Bellamy-Schedler's results for additive quiver varieties.

\subsection{Representation types and stabiliser types}\hfill

Affine GIT quotients admit a finite stratification by stabiliser types; see, for example, \cite[\S 5.2]{ShuBFM}. Beware that the general construction only gives a partition into locally closed subsets, and the smoothness of a stratum or whether its closure is a union of strata is subject to further verification. Since a general GIT quotient is defined by glueing affine GIT quotients, every multiplicative quiver variety thus admits a partition into locally closed subsets
$$
\mathcal{M}(\bs\zeta,\mathbf{d},\bs\theta)=\bigsqcup_{(H)}\mathcal{M}(\bs\zeta,\mathbf{d},\bs\theta)_{(H)},
$$
where $H$ denotes a reductive subgroup of $G=\prod_{v\in Q_0}\GL_{d_v}$, and $(H)$ denotes the conjugacy class of subgroups represented by $H$. We would like to compare this partition with another one defined by quiver representations.

Let $\mathbf{V}$ be a $\mathbf{d}$-dimensional $\bs\theta$-polystable representation. If $\mathbf{V}\cong \mathbf{V}_1^{\oplus m_1}\oplus\cdots \mathbf{V}_t^{\oplus m_t}$ is the decomposition into mutually nonisomorphic representations and each $\mathbf{V}_s$ is a $\mathbf{d}_s$-dimensional $\bs\theta$-stable representation, then we say that $\mathbf{V}$ has representation type $\bs\tau:=(m_1,\mathbf{d}_1;\ldots;m_t,\mathbf{d}_t)$.
According to Crawley-Boevey-Shaw and Schedler-Tirelli, the multiplicative quiver variety $\mathcal{M}(\bs\zeta,\mathbf{d},\bs\theta)$ admits a partition according to representation types; see \cite[\S 7]{CBS} and \cite[\S 2.7]{ST}. More precisely, we have 
$$
\mathcal{M}(\bs\zeta,\mathbf{d},\bs\theta)=\bigsqcup_{\bs\tau}\mathcal{M}(\bs\zeta,\bs\tau,\bs\theta),
$$
where $\mathcal{M}(\bs\zeta,\bs\tau,\bs\theta)$ is the locally closed subset consisting of $\bs\theta$-polystable representations of type $\bs\tau$.

\begin{Prop}\label{Prop-rep-type=stab-type}
For any $\bs\tau$ such that $\mathcal{M}(\bs\zeta,\bs\tau,\bs\theta)$ is nonempty, there is a reductive subgroup $H_{\bs\tau}\subset G$ such that the stabiliser of every $\bs\theta$-polystable representation of type $\bs\tau$ is conjugate to $H_{\bs\tau}$, and if $\bs\tau_1\ne\bs\tau_2$, then $(H_{\bs\tau_1})\ne(H_{\bs\tau_2})$. Consequently, we have 
$$
\mathcal{M}(\bs\zeta,\bs\tau,\bs\theta)=\mathcal{M}(\bs\zeta,\mathbf{d},\bs\theta)_{(H_{\bs\tau})}.
$$
\end{Prop}
\begin{proof}
The stabilisers can be explicitly described as follows. Let 
\begin{equation}\label{eq-Prop-rep-type=stab-type}
(\mathbf{V},\bs\phi)=(\mathbf{V}_1,\bs\phi_1)^{\oplus m_1}\oplus\cdots\oplus(\mathbf{V}_t,\bs\phi_t)^{m_t}
\end{equation}
represent a $\bs\theta$-polystable representation of type $\bs\tau=(m_1,\mathbf{d}_1;\ldots;m_t,\mathbf{d}_t)$ with $\dim\mathbf{V}_s=\mathbf{d}_s$ for any $s$. Then, we have
$$
\End(\mathbf{V},\bs\phi)=\gl_{m_1}\oplus\cdots\oplus\gl_{m_t};
$$
the entry of $\gl_{m_s}$ at the $x$-th row and $y$-th column represents a scalar homomorphism from the $y$-th $(\mathbf{V}_s,\bs\phi_s)$ to the $x$-th $(\mathbf{V}_s,\bs\phi_s)$. Since an element of $G$ defines an automorphism of $(\mathbf{V},\bs\phi)$ precisely when it stabilises this representation, the stabiliser of $(\mathbf{V},\bs\phi)$ is $\prod_{s=1}^t\GL_{m_s}$, which consists of the invertible endomorphisms in $\End(\mathbf{V},\bs\phi)$. Note that the conjugacy class of the stabiliser group is constant on a representation type stratum, since any $\bs\theta$-polystable representation of type $\bs\tau$ admits a conjugate of the form (\ref{eq-Prop-rep-type=stab-type}); thus, we may take $H_{\bs\tau}$ to be this subgroup. 

Let us show that if $H_{\bs\tau_1}=H_{\bs\tau_2}=:H$, then $\bs\tau_1=\bs\tau_2=:\bs\tau$; this amounts to recovering the type $\bs\tau$ from the group $H_{\bs\tau}$ above.  Let $T\subset H$ be a maximal torus which we write in the form $(\mathbb{C}^{\ast})^{m_1}\times\cdots(\mathbb{C}^{\ast})^{m_t}$. In the explicit form of $H_{\bs\tau}$ as above, a direct factor of $(\mathbb{C}^{\ast})^{m_s}$ acts by multiplication on a graded subspace $\mathbf{W}\subset\mathbf{V}$ with $\dim\mathbf{W}=\mathbf{d}_s$. This defines a grading on $\mathbf{V}$. Let us show that $\mathbf{W}$ is a subrepresentation. Indeed, for any $h=(h_v)_{v\in Q_0}\in H$ and any $a\in Q_1$, we have the following commutative diagram
$$
\begin{tikzcd}[row sep=2.5em, column sep=2.5em]
V_{t(a)} \arrow[d, "h_{t(a)}"'] \arrow[r, "\phi_a"] & V_{h(a)} \arrow[d, "h_{h(a)}"] \\
V_{t(a)}  \arrow[r, "\phi_a"] & V_{h(a)}. 
\end{tikzcd}
$$ 
Suppose that $h$ belongs to the $\mathbb{C}^{\ast}$-factor of $T$ defining $\mathbf{W}$, so that its action on $W_{t(a)}$ is the multiplication by some $\lambda\in\mathbb{C}^{\ast}$. The above diagram then implies that the action of $h_{h(a)}$ on $\phi_a(W_{t(a)})$ is also $\lambda$. We see that $\phi_a(W_{t(a)})\subset W_{h(a)}$; that is, $\mathbf{W}$ is a subrepresentation. If $\mathbf{W}$ further decomposes into a direct sum of smaller representations, then $H$, which is the stabiliser, contains a torus of a dimension larger than $\sum_sm_s$, which is not possible. It follows that each $\mathbf{W}$ as above is a stable representation. Now, the off-diagonal entries of the $\gl_{m_s}$-component of $\End(\mathbf{V},\bs\phi)$ identifies two stable subrepresentations of the form $\mathbf{W}$; therefore, the multiplicities in the type $\bs\tau$ can also be recovered from $H_{\bs\tau}$.
\end{proof}
\begin{Rem}
It is possible that for some $\bs\tau_1\ne\bs\tau_2$, we have $H_{\bs\tau_1}\cong H_{\bs\tau_2}$ but $(H_{\bs\tau_1})\ne (H_{\bs\tau_2})$, since the isomorphism type of $H_{\bs\tau}$ only depends on the multiplicities $(m_1,\ldots,m_t)$. The point of the above proof is that $(\mathbf{d}_1,\ldots,\mathbf{d}_t)$ can be recovered from $(H_{\bs\tau})$.
\end{Rem}

For any type $\bs\tau$, we write 
$$
\mathcal{M}_B(\bs\zeta,\bs\tau,\bs\theta)=\mathcal{M}(\bs\zeta,\bs\tau,\bs\theta)\cap\mathcal{M}_B(\bs\zeta,\mathbf{d},\bs\theta),
$$
so that these locally closed subvarieties form a partition of $\mathcal{M}_B(\bs\zeta,\mathbf{d},\bs\theta)$; these are the representation type strata of character varieties. Note that if $\mathbf{d}_s$ is the dimension vector of a $\bs\theta$-stable representation that is a direct summand of $(\mathbf{V},\bs\phi)\in \mathcal{M}_B(\bs\zeta,\mathbf{d},\bs\theta)$, then up to reflections the components of $\mathbf{d}_s$ are decreasing along every leg; therefore, $\mathcal{M}_B(\bs\zeta,\mathbf{d}_s,\bs\theta)$ is isomorphic to a tame character variety, and thus is connected.

\begin{Prop}
The locally closed subset $\mathcal{M}_B(\bs\zeta,\bs\tau,\bs\theta)$ is connected for every $\bs\tau$.
\end{Prop}
\begin{proof}
Suppose $g=0$, so that $\mathcal{M}_B(\bs\zeta,\mathbf{d},\bs\theta)=\mathcal{M}(\bs\zeta,\mathbf{d},\bs\theta)$. Consider the direct sum map
\begingroup
\allowdisplaybreaks
\begin{align*}
\prod_{s=1}^t\mathcal{M}(\bs\zeta,\mathbf{d}_s,\bs\theta)&\longrightarrow\mathcal{M}(\bs\zeta,\mathbf{d},\bs\theta),\\
((\mathbf{V}_s,\bs\phi_s))_{1\le s\le t}&\longmapsto \bigoplus_{s=1}^t(\mathbf{V}_s,\bs\phi_s)^{\oplus m_s}.
\end{align*}
\endgroup
The image of the stable locus $\prod_{s=1}^t\mathcal{M}(\bs\zeta,\mathbf{d}_s,\bs\theta)^s$ is precisely $\mathcal{M}(\bs\zeta,\bs\tau,\bs\theta)$. Since the domain is irreducible, its stable locus is connected, and so is its image. For $g>0$, the argument is the same; however, we need to notice that the invertibility of the $\star$-components of each $\bs\phi_s$ is equivalent to the invertibility of the $\star$-component of their direct sum.
\end{proof}
 
\subsection{Symplectic leaves of tame character varieties}\label{subsec-symp-l-char}\hfill

Let $M$ be a symplectic variety and let $\{-,-\}$ denote the Poisson bracket. A locally closed subvariety $V\subset M$ is called a symplectic leaf if it is a maximal connected Poisson subvariety on which the induced Poisson bracket is nondegenerate. By \cite[Theorem 2.3]{Kal}, there are only finitely many symplectic leaves of $M$. There is another stratification associated to a symplectic variety, the \textit{singularity stratification}; for the sequence of subvarieties $M=M_0\supset M_1\supset M_2\supset\cdots\supset M_k$ such that each $M_i$ is the singular locus of $M_{i-1}$ and $M_k$ is regular, the strata are the connected components of $M_{i-1}\setminus M_i$. By \cite[Lemma 2.5]{MGN}, this is indeed a stratification.

\begin{Thm}\label{Thm-rep-str=leaves}
The locally closed subsets $\mathcal{M}_B(\bs\zeta,\bs\tau,\bs\theta)$ indexed by representation types are precisely the symplectic leaves of $\mathcal{M}_B(\bs\zeta,\mathbf{d},\bs\theta)$. In particular, 
\begin{itemize}
\item each $\mathcal{M}_B(\bs\zeta,\bs\tau,\bs\theta)$ is smooth, and
\item the disjoint union $\mathcal{M}_B(\bs\zeta,\mathbf{d},\bs\theta)=\bigsqcup_{\bs\tau}\mathcal{M}_B(\bs\zeta,\bs\tau,\bs\theta)$ is indeed a stratification.
\end{itemize}
\end{Thm}
\begin{proof}
Write $\mathcal{M}:=\mathcal{M}_B(\bs\zeta,\mathbf{d},\bs\theta)$. Consider the formal neighbourhood $f_{\rho}:\widehat{\mathcal{M}}_{\rho}\rightarrow\mathcal{M}$ for any $\rho\in \mathcal{M}$. We will show that the inverse images under $f_{\rho}$ of the representation type strata are precisely the singularity strata, so that the two stratifications coincide. The theorem then follows from \cite[Lemma 2.5]{MGN}. Identify $\widehat{\mathcal{M}}_{\rho}$ with the formal neighbourhood of the zero representation $0$ in an additive quiver variety $\mathfrak{M}$, denoted by $\widehat{\mathfrak{M}}_0$. By \cite[Theorem 5.16 (2)]{KS}, the inverse image of a stabiliser type stratum is a strabiliser type stratum under the natural morphism $\widehat{\mathfrak{M}}_0\rightarrow \mathcal{M}$. However, a stabiliser type stratum in $\mathfrak{M}$ is a representation type stratum by Proposition \ref{Prop-rep-type=stab-type}. By \cite[Proposition 3.15]{BS0}, the stabiliser type strata of $\widehat{\mathfrak{M}}_0$ are the symplectic leaves and thus are the singularity strata. The theorem is proved.
\end{proof}

The classification of codimension-two symplectic leaves can be achieved exactly as in \cite[Section 7]{BS0}. The formal neighbourhood of any $\rho\in\mathcal{M}(\bs\zeta,\mathbf{d},\bs\theta)$ can be described by an additive quiver variety associated to the ext-quiver of $\rho=(\mathbf{V},\bs\phi)$, which we now recall. Suppose that $(\mathbf{V},\bs\phi)$ is of type $\bs\tau$, and we write $(\mathbf{V},\bs\phi)\cong\bigoplus_{s\in Q_0'}(\mathbf{V}_s,\bs\phi_s)^{m_s}$ for some finite set $Q_0'$ with $(\mathbf{V}_s,\bs\phi_s)\ncong(\mathbf{V}_{s'},\bs\phi_{s'})$ for $s\ne s'$ with each $(\mathbf{V}_s,\bs\phi_s)$ being a $\bs\theta$-stable representation. Denote by $\mathbf{d}_s$ the dimension vector of $\mathbf{V}_s$. The double of the ext quiver has $Q'_0$ as its set of vertices and $-(\mathbf{d}_{s_1},\mathbf{d}_{s_2})$ arrows from $s_1$ to $s_2\ne s_1$, and will be denoted by $Q'$. Denote by $\mathfrak{M}$ the additive quiver variety defined by the quiver $Q'$ and dimension vector $\mathbf{m}=(m_s)_{s\in Q_0'}$ with trivial deformation parameter and trivial stability condtion. The formal neighbourhood $\smash{\widehat{\mathfrak{M}}_0}$ of the zero representation in $\mathfrak{M}$ is then isomorphic to $\smash{\widehat{\mathcal{M}}_{\rho}}$ according to \cite[Theorem 5.11]{Dav24} or \cite[Theorem 5.4]{KS}.

We say that $\bs\tau$ is an isotropic type if the following conditions hold:
\begin{itemize}
\item If $\mathbf{d}_s$ is an imaginary root, then $m_s=1$.
\item If $Q''$ is the quiver obtained by removing the loops from the Ext-quiver $Q'$, then $Q''$ is an affine Dynkin quiver.
\item The dimension vector $\mathbf{m}$ is the minimal positive imaginary root.
\end{itemize}

\begin{Prop}\label{Prop-codim2-leaves}
Let $\mathbf{d}\in \Sigma_{\bs\zeta,\bs\theta}$ be an imaginary root. Then, the codimension two strata of $\mathcal{M}_B(\bs\zeta,\mathbf{d},\bs\theta)$ are in bijection with the isotropic types. For any $\rho\in\mathcal{M}_B(\bs\zeta,\bs\tau,\bs\theta)$ with $\bs\tau$ being isotropic, the formal transversal slice at $\rho$ is isomorphic to the formal neighbourhood of zero in the additive quiver variety $\mathfrak{M}$ associated to $Q''$.
\end{Prop}
\begin{proof}
Bellamy-Schedler's proof for additive quiver varieties remains valid in the multiplicative case. The transversal slice computations only involve the additive quiver variety, and thus  \cite[Lemma 7.1 and Proof of Theorem 1.20]{BS0} directly apply. The proof of \cite[Lemma 7.2]{BS0} relies on the dimension formula for the representation type strata, which has the same form in the multiplicative case (see \cite[Proposition 2.15]{ST}).
\end{proof}

\addtocontents{toc}{\protect\setcounter{tocdepth}{-1}}
\bibliographystyle{alpha}
\bibliography{BIB}
\end{document}